\documentclass[11pt,reqno]{amsart}
\usepackage{cite}
\usepackage{amssymb}
\usepackage{amsfonts}

\newcommand{\rmv}[1]{}
\def\dsd{\mathrel{\ast\ast}}
\newcommand{\set}[1]{\left\{#1\right\}}
\def\FF{\mathbb{F}}
\def\GGM{\mathsf{GGM}}
\def\Radf#1#2{\mathrm{Rad}_{#1}(#2)}

\def\MM{\ensuremath{{\mathcal M}}}
\def\pv#1{\ensuremath{{\bf#1}}}

\def\QQ{\ensuremath{\mathbb Q}}
\def\FF{\ensuremath{\mathbb F}}

\def\onto{\twoheadrightarrow}

\def\inv{^{-1}}
\def\p{\varphi}

\def\pinv{{\p \inv}}
\def\J{\mathrel{{\mathcal J}}} % J - relation
\def\R{\mathrel{{\mathcal R}}} % R - relation
\def\L{\mathrel{{\mathcal L}}} % L - relation
\def\H{\mathrel{{\mathcal H}}} % H - relation

\def\e<{\leq _{E}}
\def\j{\mathrel{\leq _{\J}}}

\def\ov#1{\ensuremath{\overline {#1}}}

\def\til#1{\ensuremath{\widetilde {#1}}}

\def\malce{\mathop{\hbox{$\bigcirc$\kern-8.5pt\raise1pt
\hbox{\scriptsize$m$}\kern1.5pt}}}
\def\1sk{^{(1)}}

\def\to{\rightarrow}

\def\End#1#2{\mathrm{End}_{#1} (#2)}
\def\apex#1{\mathrm{Apx}(#1)}

\def\Thmname{Theorem}
\def\Propname{Proposition}
\def\Lemmaname{Lemma}
\def\Definitionname{Definition}
\newtheorem{Thm}{\Thmname}[section]
\newtheorem{Prop}[Thm]{\Propname}
\newtheorem{Lemma}[Thm]{\Lemmaname}

\newtheorem{Cor}[Thm]{Corollary}

\newtheorem{Claim}{Claim}

\numberwithin{equation}{section}

\title[Representation Theory of Finite Semigroups]
{Representation Theory of Finite Semigroups, Semigroup
Radicals and Formal Language Theory}

\author[J. Almeida\and S. Margolis\and B. Steinberg\and
M. Volkov]{Jorge Almeida\and Stuart Margolis\and Benjamin
Steinberg\and Mikhail Volkov}
\address{J. Almeida: Departamento de Matem\'atica Pura,
Faculdade de Ci\^encias, Universidade do Porto, 4169-007 Porto,
PORTUGAL} \email{jalmeida@fc.up.pt}
\address{S. Margolis: Department of Mathematics,
Bar Ilan University, 52900 Ramat Gan, ISRAEL}
\email{margolis@math.biu.ac.il}
\address{B. Steinberg: School of Mathematics and Statistics,
Carleton University Ottawa, Ontario K1S 5B6, CANADA}
\email{bsteinbg@math.carleton.ca}
\address{M. Volkov: Department of Mathematics and Mechanics,
Ural State University, 620083 Ekaterinburg, RUSSIA}
\email{Mikhail.Volkov@usu.ru}
\thanks{The first author acknowledges the support of the Centro de
Matem\'atica da Universidade do Porto, financed by FCT through the
programmes POCTI and POSI, with Portuguese and European Community
structural funds.
  The second author acknowledges the support of the Excellency Center,
``Group Theoretic Methods for the Study of Algebraic Varieties''
of the Israeli Science Foundation and thanks Professor J.-\'E. Pin
for inviting him to be a visitor to LIAFA. The third author
acknowledges the support of NSERC. The fourth author acknowledges
support from the Russian Foundation for Basic Research, grant
05-01-00540.}

\keywords{representation theory, radicals, language theory}
\subjclass{}

\begin{document}

\begin{abstract}
In this paper we characterize the congruence associated to the
direct sum of all irreducible representations of a finite semigroup
over an arbitrary field, generalizing results of Rhodes for the
field of complex numbers. Applications are given to obtain many new
results, as well as easier proofs of several results in the
literature, involving: triangularizability of finite semigroups;
which semigroups have (split) basic semigroup algebras, two-sided
semidirect product decompositions of finite monoids; unambiguous
products of rational languages; products of rational languages with
counter; and \v{C}ern\'y's conjecture for an important class of
automata.
\end{abstract}

\maketitle

{\small
\tableofcontents

}

\section{Introduction}

For over 100 years, the theory of linear representations has
played a fundamental role in studying finite groups, finite
dimensional algebras and Lie algebras as well as other parts of
algebra. By way of contrast, the theory of semigroup
representations, which was intensively developed during the 50s
and 60s in classic work such as Clifford \cite{Cliffordold},
Munn~\cite{Munnalg,Munn} and Ponizovsky (see~\cite[Chapter~5]{CP}
for an account of this work, as well as
\cite{Lallement,RhodesZalc} for nicer treatments restricting to
the case of finite semigroups), has found almost no applications
in the theory of finite semigroups. It was pointed out by
McAlister in his survey of 1971~\cite{McAlister} that the only
paper applying representation theoretic results to finite
semigroups was the paper~\cite{rhodeschar} of Rhodes. This paper
determined the congruence on a finite semigroup $S$ associated to
the direct sum of the irreducible representations of $S$ over the
field of complex numbers. Rhodes applied this result to calculate
the Krohn-Rhodes complexity~\cite{KR2,Arbib} of completely regular
semigroups. Around the time of McAlister's survey, there also
appeared a paper of Zalcstein~\cite{zalc} trying to apply
representation theory to finite semigroup theory.

For many years, representation theory of finite semigroups remained
dormant until Putcha, in a series of papers (cf.\
\cite{Putcharep,Putchaquiver, Putchaweights,PutchaRecip} and
others), revived the theme. Putcha was primarily interested in
relating semigroup theory with modern areas in representation theory
such as quasi-hereditary algebras, weights for representations of
finite groups of Lie type and with calculating quivers of algebras
of finite semigroups.  However, his research was not aimed at
applying the representation theory of semigroups to the study of
finite semigroups in their own right.   While to some extent we
continue in the vein of relating finite semigroup theory to the rest
of modern algebra --- for instance we determine over an arbitrary
field $K$ which finite semigroups have basic or split basic
semigroup algebras over $K$ --- we very much focus on using
representation theory precisely for the purpose of answering
questions from finite semigroup theory. We are particularly
interested in varieties of finite semigroups and their connections
with formal language theory and other aspects of theoretical
computer science, as exposited in the two treatises by Eilenberg
\cite{Eilenberg}.  Nonetheless we expect that the first four
sections of this paper should be of interest to readers in Algebraic
Combinatorics, Representation Theory and Finite Semigroup Theory.

Let us briefly survey the contents of the paper. Following the
preliminaries, we define the Rhodes radical of a finite semigroup
$S$ with respect to a field $K$ to be the congruence on $S$
induced by the Jacobson radical of its semigroup algebra $KS$ over
$K$. Using classical Wedderburn theory, we give a more conceptual
proof of Rhodes's characterization~\cite{rhodeschar} of this
radical in characteristic zero and extend it to characteristic
$p$. Further, the radical is shown to be intimately related with
the Mal'cev product, which is an integral part of the varietal
theory of finite semigroups.  We also give an alternative
semigroup representation theoretic proof of the description of the
Rhodes radical, along the original lines of Rhodes
\cite{rhodeschar}, that allows a more precise and usable
characterization of the radical.

Part of our aim is to render things in a form understandable to
both specialists and non-specialists.  Recent work of Bidigaire
\textit{et al.} \cite{BHR} and Brown \cite{Brown1,Brown2}, for
instance, spends quite some time in redeveloping basic aspects of
the representation theory of idempotent semigroups (known as
bands) that were already in the literature
\cite{McAlister,RhodesZalc,MunnIdem,Lallement}, but perhaps not in
a form accessible to most mathematicians.  Our results
handle the general case in a form that both semigroup theorists as
well as workers in finite dimensional algebras, group
representation theory and other related fields should find useful.

We then proceed to applications.  The first application gives
abstract, algebraic characterizations of finite semigroups that are
triangularizable over a field $K$; in the language of the theory of
finite dimensional algebras, we characterize those finite semigroups
whose semigroup algebras are split basic $K$-algebras. The case of a
finite field was handled by three of the authors in~\cite{amv}
without using representation theory, leading to a much more
complicated proof. Here we handle all fields $K$ in a uniform manner
by simply characterizing those semigroups all of whose irreducible
representations over $K$ have degree one. It turns out that the
collection of finite semigroups triangularizable over a given field
$K$ is a variety of finite semigroups (that of course depends on
$K$) and that those ``triangularizable'' varieties are in fact some
of the most commonly studied varieties in finite semigroup theory.

Our next application is to obtain simpler proofs of some bilateral
semidirect product decomposition results of Rhodes, Tilson and
Weil~\cite{Kernel,MPS2} using representation theory.  The original
proofs rely on a case-by-case analysis of Rhodes's classification
of maximal proper surjective morphisms.

After purely algebraic applications, we switch to those dealing
with important objects of theoretical computer science such as
formal languages and finite automata. We use modular
representation theory to give simpler proofs of results of
P\'eladeau and Weil~\cite{collapse,Weil} on marked products with
modular counter and characteristic zero representation theory to
obtain simpler proofs of results of Pin, Straubing and
Therien~\cite{Pinetal} on unambiguous marked products. Our final
application uses representation theory to confirm the longstanding
\v{C}ern\'y conjecture on synchronizing automata in the special
case that the transition monoid belongs to the much-studied
variety \pv{DS}.

Further applications of our results have been obtained by the third
author \cite{Mobius,MobiusII}; in particular the results of Bidigaire
\emph{et al.} \cite{BHR} and Brown \cite{Brown1,Brown2} on random
walks on hyperplane arrangements and on bands have been extended to
the varieties \pv {DA} and $\pv {DO}\cap \ov {\pv
{Ab}}$, which is as far as these results can be extended.

We have tried to make the representation part (Section~3) of this
paper accessible both to readers from semigroup theory and readers
familiar with representation theory from other contexts. Having
the latter category of readers in mind, in the next section we
give a concise overview of standard notions and terminology of
semigroup theory needed for the representation part. The
application part (Sections~4--7) requires further background in
semigroup theory, formal languages and automata.

\section{Preliminaries}
Good sources for semigroup theory, in particular finite semigroup
theory, are~\cite{CP,Arbib,Eilenberg,Pin,almbook}. Here we introduce
some standard notions and terminology. The reader is welcome to skip
this section, referring back only as needed.

A \emph{congruence} on a semigroup $S$ is an equivalence
relation $\equiv$ such that
$$s\equiv s'\implies ts\equiv ts',\ st\equiv s't$$
for all $s,s',t\in S$. Left and right congruences are defined
analogously. If $\p:S\to T$ is a morphism, then the congruence
\emph{associated to $\p$} is defined by $s\equiv_\p t$ if and only
if $s\p=t\p$.

An \emph{idempotent} $e$ of a semigroup is an element such that
$e=e^2$. It is well known that in a finite semigroup $S$ some power
of each element is an idempotent; namely, for all $s\in S$, one
verifies that $s^{|S|!}$ is idempotent. The set of idempotents of a
semigroup $S$ is denoted $E(S)$. It is a partially ordered set via
the order
\begin{equation}\label{order}
e\leq f\iff ef=fe=e.
\end{equation}
A \emph{semilattice} $E$ is an idempotent commutative semigroup. In
this case, the order \eqref{order} has all finite meets, the meet
being given by the product in $E$.

A \emph{right ideal} of a semigroup $S$ is a subset $R$ such that
$RS\subseteq R$.  Left ideals and (two-sided) ideals are defined
similarly.  If $s\in S$, we use $R(s)$, $L(s)$, $J(s)$ for the
respective right, left and two-sided principal ideals generated by
$s$. This leads to the definitions of Green's relations
\cite{Green,Arbib,CP}, which play an essential role in semigroup
theory. We define an equivalence relation $\R$ on $S$ by setting,
for $s,t\in S$,  $s\R t$ if and only if $R(s)=R(t)$; in this case
one writes $R_s$ for the $\R$-class of $s$. One similarly defines
the equivalence relations $\L$ and $\J$, whose classes of $s$ are
denoted $L_s$ and $J_s$ respectively.  Define $s\H t$ if $s\R t$
and $s\L t$; the $\H$-class of $s$ is denoted $H_s$.  There are
also associated preorders. For instance, $s\leq_{\R} t$ if and
only if $R(s)\subseteq R(t)$. It is easy to see that $\R$ is a
left congruence and $\L$ is a right congruence.

In a finite semigroup (or even in an algebraic semigroup
\cite{PutchaAlg,Renner}), the following \emph{stability} relations
hold~\cite{Arbib}:
$$s\J st\iff s\R st,\quad t\J st\iff t\L st.$$  From these relations,
it follows that in a finite semigroup, if $s\J t$ then there exists
$u\in S$ such that $s\R u\L t$ and $v\in S$ such that $s\L v\R t$.  In
the case that $J_s$ is a subsemigroup one can take $u$ and $v$ to be
idempotents as every $\H$-class within $J_s$ contains an idempotent.

An element $s\in S$ is called (\emph{von Neumann}) \emph{regular} if
$s\in sSs$.  In a finite semigroup, $s$ is regular if and only if $J_s$
contains an idempotent if and only if $R_s$ contains an idempotent if
and only if $L_s$ contains an idempotent.  A $\J$-class (respectively,
$\R$-class, $\L$-class) is called \emph{regular} if it contains an
idempotent. If $e$ is an idempotent, then $H_e$ is a group, called
the \emph{maximal subgroup} at $e$. It is the group of units of
the \emph{local monoid} $eSe$ and so it is the largest subgroup
of $S$ with identity $e$.  By a subgroup of a semigroup $S$, we
mean simply a subsemigroup that is a group; it need not have
the same identity as $S$ in the case that $S$ is a monoid.  The local
monoid $eSe$ is the largest subsemigroup of $S$ with identity $e$.
For example, if $S$ is the monoid of $n\times n$ matrices over $K$ and
$e$ is an idempotent of rank $r$, then $eSe$ is isomorphic to the
monoid of $r\times r$ matrices over $K$ and the maximal subgroup $H_e$
is isomorphic to the general linear group of degree $r$ over $K$.

If $S$ is a semigroup, we set $S^1$ to be $S$ with an adjoined
identity if $S$ is not a monoid and $S$ otherwise.   We shall
frequently use the following fact: suppose that $e,f\in E(S)$; then
\begin{equation}\label{leftrightidentity}
\begin{split}
e\L f\Leftrightarrow ef = e,\ fe =f\\
e\R f\Leftrightarrow ef = f,\ fe =e
\end{split}
\end{equation}
For instance, if $e\L f$, then $e= xf$ for some $x\in S^1$. Hence
$$ef = xff = xf =e;$$ the other equalities are handled similarly.

A semigroup is called \emph{simple} if it has no proper
(two-sided) ideals. A semigroup $S$ with $0$ is called
$0$-\emph{simple} if $S^2=\{st\mid s,t\in S\}\ne0$ and the only
ideals of $S$ are $\{0\}$ and $S$. Simple semigroups and
$0$-simple semigroups were classified up to isomorphism by Rees
and Suschewitsch~\cite{CP}. We shall need in the sequel only the
following properties that are the content of~\cite[XI.
Propositions 1.2--1.4]{TilsonEilen}.

\begin{Prop}\label{simplesemigroups}
Let $S$ be a finite simple semigroup.  Then every element of $S$
belongs to a subgroup of $S$. For any idempotents $e,f\in S$,
there exist $x\in eSf$ and $y\in fSe$ such that $e=xy$ and $f=yx$.
Moreover, $eSe$ is the group $H_e$, $fSf$ is the group $H_f$ and
the map $H_e\to H_f$ given by $h\mapsto yhx$ is a group
isomorphism.
\end{Prop}

This proposition says that the idempotents of a simple semigroup are
\emph{conjugate} and that the local monoids are the maximal subgroups;
moreover, they are all isomorphic to the same group.

An ideal of a semigroup $S$ is called \emph{minimal} if it
contains no other ideal of $S$; the minimal ideal of a finite
semigroup is a simple semigroup and is a regular $\J$-class
~\cite{CP}. An ideal of a semigroup $S$ with $0$ is called
\emph{$0$-minimal} if the only ideal of $S$ properly contained in
it is $\{0\}$; a $0$-minimal ideal $I$ of a finite semigroup is
either $0$-simple (and then $I\setminus\{0\}$ is a regular
$\J$-class) or it is \emph{null}, meaning $I^2=0$.

The following definition, introduced by Eilenberg and
Sch\"utzenberger~\cite{Eilenschutz,Eilenberg}, is crucial in
finite semigroup theory. A class \pv V of semigroups closed under
formation of finite direct products, subsemigroups and homomorphic
images is called a \emph{variety of  finite semigroups} (or
sometimes a \emph{pseudovariety} of semigroups). Varieties of
finite monoids and groups are defined analogously.  We remark that
in universal algebra, the term variety is used differently, but
since we shall not consider such varieties, no confusion should
arise.

Some varieties that shall play an important role in this paper are
the trivial variety \pv I (containing only the trivial semigroup) and
the variety of finite $p$-groups ($p$ a prime) $\pv {G}_p$.  The
variety of finite Abelian groups is denoted \pv {Ab}. The variety
of finite semilattices is denoted $\pv{Sl}$.

The following notion shall be used throughout this paper. If \pv V
is a variety of finite semigroups, a morphism $\p:S\to T$ is called
a \pv V-\emph{morphism} if, for each idempotent $e\in T$, its
preimage $e\pinv$ (which is then a subsemigroup in $S$) belongs to
\pv V. The congruence associated to a \pv V-morphism is called a \pv
V-\emph{congruence}. In other words, a congruence on $S$ is a \pv
V-congruence if and only if all its congruence classes that are
subsemigroups belong to \pv V.  For instance, if $\p:G\to H$ is a
group homomorphism, then $\p$ is a \pv V-morphism if and only if
$\ker \p\in \pv V$.

Finally, we recall two fundamental varietal constructions. If \pv
V is a variety of finite monoids, then \pv {LV} denotes the class
of all finite semigroups $S$ such that, for each idempotent $e\in
S$, the local monoid $eSe$ belongs to \pv V. It is easy to see
that \pv {LV} is a variety of finite semigroups. If \pv V and \pv
W are varieties of finite semigroups, their \emph{Mal'cev product}
$\pv V\malce \pv W$ consists of all finite semigroups $S$ such
that there is a finite semigroup $T$ mapping homomorphically onto $S$
such that
$T$ admits a $\pv V$-morphism to a semigroup in $\pv W$. Again, it
is well known and easy to verify that $\pv V\malce \pv W$ is a
variety of finite semigroups.

The most important example is when \pv V is a variety of finite
groups.  A semigroup is a \emph{local group} if $eSe$ is a group
for each idempotent $e$.  For instance, by
Proposition~\ref{simplesemigroups} simple semigroups are local
groups. If one considers all $n\times n$ upper triangular matrices
over a field $K$ that have a fixed zero/non-zero pattern on the
diagonal, we will see in Section 4 that one obtains a local group.
Thus the monoid of all upper triangular matrices is a disjoint
union of local groups.

Our goal is to state the well known version of
Proposition~\ref{simplesemigroups} for local groups.
Unfortunately, even though this is folklore in semigroup theory,
we could not pinpoint an exact reference.  First we need the
following well-known finiteness result, which is a ``Pumping
Lemma'' for finite semigroups~\cite[Proposition 5.4.1]{almbook}.
Set $S^n$ be the ideal of $S$ consisting of all elements of $S$
that can be expressed as a product of $n$ elements of $S$.

\begin{Lemma}[Pumping Lemma]\label{pumpinglemma}
Let $S$ be a semigroup with $n$ elements.  Then $S^n=SE(S)S$.
\end{Lemma}

Now we can state the main property of local groups.
\begin{Prop}\label{localgroups}
Let $S$ be a finite semigroup.  Then the following are equivalent:
\begin{enumerate}
\item $S$ is a local group;
\item $S^n$ is a simple semigroup for some $n>0$ (i.e.\ $S$ is a
  nilpotent ideal extension of a simple semigroup);
\item $S^n$ is the minimal ideal of $S$ for some $n>0$ (i.e.\ $S$ is a nilpotent
  extension of its minimal ideal);
\item $S$ does not contain a semigroup isomorphic to the two-element
  semilattice $\{0,1\}$ with
  multiplication.
\end{enumerate}
Furthermore for any idempotents $e,f\in S$, there exist $x\in eSf$
and $y\in fSe$ such that $e=xy$ and $f=yx$.  Moreover, the groups
$eSe$ and $fSf$ are isomorphic via the map $eSe\to fSf$ given by
$h\mapsto yhx$.
\end{Prop}
\begin{proof}
Suppose first that (1) holds.  Then (4) must hold since if
$\{e,f\}\subseteq S$ is isomorphic to $\{0,1\}$ with $e$ as the
identity, then $efe=f$ and so $e,f\in eSe$, showing that $eSe$ is not
a group.

 For (4) implies (3), let $I$ be the minimal ideal of $S$.
We show that $E(S)\subseteq I$.  Suppose $e\in E(S)\setminus I$.  Let
$s\in I$ be any element and set $n=|S|$.  Then $f=(ese)^{n!}$ is an idempotent
belonging to $I$ (so in particular $f\neq e$) and $ef=fe=f$.  Thus
$\{e,f\}$ is a subsemigroup isomorphic to $\{0,1\}$.  This
contradiction shows that $E(S)\subseteq I$.  Now by the Pumping Lemma,
if $n=|S|$, then $S^n= SE(S)S\subseteq I$.  However, $SE(S)S$ is
clearly an ideal, so $I\subseteq SE(S)S$.  Hence $S^n=I$.

We noted that the minimal ideal of any finite semigroup is a
simple semigroup so the implication (3) implies (2) is trivial.
For (2) implies (1), suppose that $T=S^n$ is a simple semigroup.
Notice that $E(S)\subseteq T$ and that $eSe\subseteq T$ for any
idempotent $e$ since $e\in S^n$ for all $n$.  Thus
\[eSe=e(eSe)e\subseteq eTe\subseteq eSe,\] so $eSe=eTe$.  But $eTe$ is
a group by Proposition~\ref{simplesemigroups}.  This proves (1).

The proof that (2) implies (1) shows that in a local group $S$ with
minimal ideal $I$, one has that $I$ contains all the idempotents of
$S$ and $eSe=eIe$ for each idempotent $e$ of $S$.
Proposition~\ref{simplesemigroups} then implies the final statement of
the proposition.
\end{proof}

\section{The Rhodes Radical}

\subsection{Background and Motivation}\label{motivation}
Let $K$ be a field and $S$ a finite semigroup. The semigroup
algebra of $S$ over $K$ is denoted $KS$. Recall that this is the
$K$-vector space with basis $S$ and the multiplication extending
the multiplication in $S$. If $A$ is a finite dimensional
$K$-algebra (for instance $KS$), then it has a largest nilpotent
ideal $\mathrm {Rad}(A)$, called its (Jacobson) \emph{radical}.
Consider the composite mapping
$$S\to KS\to KS/\mathrm{Rad}(KS);$$
this is a morphism of semigroups where the latter two are viewed
with respect to their multiplicative structure. We define $\Radf K
S$, called the \emph{Rhodes radical} of $S$, to be the associated
congruence on $S$. Let us briefly discuss the role of the Rhodes
radical for the representation theory of finite semigroups.

Let $V$ be a $K$-vector space of finite dimension $n$. Then $\End K
V$ denotes the monoid of $K$-endomorphisms of $V$. We shall identify
$\End K V$ with the monoid $M_n(K)$ of $n\times n$ matrices over $K$
whenever it is convenient. A \emph{representation} of a finite
semigroup $S$ over $K$ of \emph{degree} $n$ is a morphism
$\rho:S\to M_n(K)$ or, equivalently, a morphism $\rho:S\to \End K V$
where $V$ is an $n$-dimensional vector space over $K$. It is easy to
see that via $\rho$ we can view $V$ as a finite dimensional (right)
$KS$-module and that all finite dimensional (right) $KS$-modules
arise in this way. The \emph{regular representation} of $S$ is the
faithful representation on the $K$-vector space with basis $S^1$ and
where the action is induced by right multiplication of $S$ on the
basis elements.

A subsemigroup $S$ of $\End K V$ is called \emph{irreducible} if
there is no proper, non-zero subspace of $V$ that is invariant under
$S$. A representation $\rho:S\to \End K V$ of a semigroup $S$ is
called \emph{irreducible} if $S\rho$ is an irreducible subsemigroup
of $\End K V$.  A representation is irreducible if and only if the
associated $KS$-module is simple.

It is well known that the radical $\mathrm{Rad}(A)$ of a finite
dimensional $K$-algebra $A$ is the intersection of the kernels of
the irreducible representations of $A$. Since every irreducible
representation $\rho:S\to \End K V$ of a finite semigroup $S$
uniquely extends to an irreducible representation of the semigroup
algebra $KS$, and vice versa, every irreducible representation of
$KS$ restricts to an irreducible representation of $S$, we
conclude that the Rhodes radical $\Radf K S$ of $S$ is precisely
the intersection of the congruences of the form $\equiv_\rho$
where $\rho:S\to \End K V$ is an irreducible representation of
$S$. Thus, the Rhodes radical in the finite semigroup setting
naturally corresponds to the (Jacobson) radical in the setting of
finite dimensional algebras. Moreover, in spite of the fact that
the irreducible representations of $S$ and $KS$ are basically the
same objects, we will see that working with the Rhodes radical
$\Radf K S$ has some advantages over considering the radical
$\mathrm{Rad}(KS)$ of the corresponding semigroup algebra. The
point is that, as we are going to show, the Rhodes radical $\Radf
K S$ can be explicitly calculated in terms that are internal with
respect to the semigroup $S$ while determining the radical
$\mathrm{Rad}(KS)$ requires studying invertibility of certain
matrices in the matrix ring over the algebras $KH$ for all maximal
subgroups $H$ of $S$ (cf.\ \cite{hall}) which is, generally
speaking, a highly non-trivial task.

Rhodes~\cite{rhodeschar} calculated $\Radf K S$ for $K$ the field of
complex numbers, but his arguments work for any field of
characteristic $0$.  Extensions of these results in a more general
context have been obtained by Okni\'nski~\cite{Okninski}, but
without the varietal viewpoint~\cite{Eilenberg} that we use to tie
the results to language theory. Here we furnish two descriptions
of the Rhodes radical. The first proceeds via an argument using the
theory of finite dimensional algebras. Afterwards we give a
description along the lines of Rhodes~\cite{rhodeschar}, using
semigroup representation theory developed by Clifford, Munn and
Ponizovsky~\cite{CP,RhodesZalc,Lallement} and the semi-local theory
of Krohn, Rhodes and Tilson~\cite{Arbib}. Both proofs are
informative, the first being technically easier, the second giving a
more concrete description of the congruence.

Given a field $K$, let
$$\pv G_{K}=
\begin{cases}
  \pv I & char\ K =0\\
  \pv {G}_p & char\ K =p.
\end{cases}$$
It is well known that this is the variety of finite groups that are ``unipotent'' over $K$
(i.\,e.\  a finite group $G$ has a faithful unitriangular representation over $K$
if and only if $G\in \pv G_K$). This notation will allow us to phrase
our results in a characteristic-free manner.

We shall also often encounter the variety $\pv{LG}_{K}$.  By
Proposition~\ref{localgroups} a finite semigroup $S$ belongs to this
variety if and only if there is an integer $n$
such that $S^n$ is a simple
semigroup $U$, all
of whose maximal subgroups are in $\pv G_{K}$. Equivalently, $S\in
\pv{LG}_{K}$ if and
only if it does not contain a copy of the two element semilattice
$\{e,f\mid ef=fe=e^2=e,\
f^2=f\}$ and the maximal subgroups of $S$ belong to $\pv
G_{K}$.

\subsection{Rhodes Radical via Wedderburn Theory}\label{augmentation}

Our first goal is to relate the notion of a $\pv V$-morphism
to algebra morphisms.

\begin{Lemma}\label{toalgebra}
Let $\p:A\to B$ be a morphism of $K$-algebras with $\ker \p$
nilpotent.  Let $S$ be a finite subsemigroup of $A$.  Then $\p|_S$
is an $\pv{LG}_{K}$-morphism.
\end{Lemma}

\begin{proof}
Without loss of generality, we may assume that $S$ spans $A$ and
hence that $A$ is finite dimensional. Let $e_0\in E(B)$ and
$U=e_0\p|_S\inv$. First we show that $U$ does not contain a copy
of the two element semilattice. Indeed, suppose that $e,f\in E(U)$
and $ef=fe=e$. Then
$$(f-e)^2 = f^2 -ef -fe+e^2 = f-e.$$
Since $f-e\in \ker \p$, a nilpotent ideal, we conclude $f-e=0$,
that is $f=e$.

Now let $G$ be a maximal subgroup of $U$ with identity $e$. Then
$g-e\in \ker \p$. Since $g$ and $e$ commute, if the characteristic
is $p$, then, for large enough $n$, $$0 = (g-e)^{p^n} =
g^{p^n}-e$$ and so $G$ is a $p$-group.  If the characteristic is
$0$, then we observe that \mbox{$(g-e)^n=0$} for some $n$ (take
$n$ minimal).  So by taking the regular representation $\rho$ of
$G$, we see that $g\rho$ is a matrix with minimal polynomial of
the form $(g\rho-1)^n$; that is $g\rho$ is unipotent.  A quick
consideration of the Jordan canonical form for such $g\rho$ shows
that if $g\rho\neq 1$, then it has infinite order. It follows that
$g=e$ and so $G$ is trivial. This completes the proof that $U\in
\pv {LG}_{K}$.
\end{proof}

Let $\p:S \to T$ be a morphism and let $\overline{\p}:KS \to KT$
denote the linear extension of $\p$ to the semigroup algebra $KS$.
Our goal is to prove the converse of Lemma~\ref{toalgebra} for
$\overline{\p}$.  Of particular importance is the case where $T$
is the trivial semigroup.  In this case $\ker \overline \p$ is
called the \emph{augmentation ideal}, denoted $\omega KS$, and
$\ov \p$ the \emph{augmentation map}. It is worth observing that
if $U$ is a subsemigroup of $S$, then the augmentation map for $U$
is the restriction of the augmentation map of $S$ and hence
$\omega KU=\omega KS\cap KU$.   So we begin by giving a
varietal characterization of finite semigroups with nilpotent
augmentation ideal.

First we prove a classical lemma showing how to find generators
for the ideal $\ker \overline{\p}$ in terms of $\p:S\to T$.
\begin{Lemma}\label{gens}
Let $\p:S \to T$ be a morphism and let $\overline{\p}:KS \to KT$
denote the linear extension of $f$ to the semigroup algebra $KS$.
Then the set
$$X = \set{s_1 - s_2\mid s_1\p=s_2\p}$$
generates the ideal $\ker \overline{\p}$ as a vector space over $K$.
\end{Lemma}
\begin{proof}
Clearly, $X \subseteq \ker \overline{\p}$. Now take an arbitrary
$u=\sum_{s\in S}c_ss\in\ker \overline{\p}$ where $c_s\in K$.
Applying the morphism $\overline{\p}$ to $u$, we obtain
$$0=\sum_{t\in S\p}(\sum_{s\p=t}c_s)t$$
whence for each $t \in S\p$,
\begin{equation}
\label{sum}
\sum_{s\p=t}c_s=0
\end{equation}
as elements of $T$ form a basis of $KT$. Now picking for each $t \in S\p$
a representative $s_t\in S$ with $s_t\p = t$ and using \eqref{sum}, we can
rewrite the element $u$ as follows:
$$u=\sum_{t\in S\p}(\sum_{s\p=t}c_s(s-s_t)),$$
that is, as a linear combination of elements in $X$.
\end{proof}

We recall a standard result from the theory of finite-dimensional
algebras due to Wedderburn~\cite{Wedderburn}.

\begin{Lemma} \label{Wed}
Let $A$ be a finite dimensional associative algebra over a field
$K$. Assume that $A$ is generated as a $K$-vector space  by a set
of nilpotent elements. Then $A$ is a nilpotent algebra.
\end{Lemma}

The following can be proved using representation theory of finite
semigroups or extracted from a general result of
Ovsyannikov~\cite{Ovs}. We give a simple direct proof using the
above lemma. A similar proof for groups can be found, for
instance, in~\cite{Passman}.

\begin{Prop} \label{simp}
Let $S$ be a finite semigroup.  Then the augmentation ideal
$\omega KS$ is nilpotent if and only if $S\in \pv {LG}_{K}$.
\end{Prop}
\begin{proof}
Recall that $\omega KS$ is the kernel of the morphism $KS
\rightarrow K$ induced by the trivial morphism $S \rightarrow
\set{1}$; it consists of all elements $\sum_{s \in S}c_{s}s \in
KS$ such that $\sum_{s \in S}c_{s}=0$.

Suppose first that $S\in \pv {LG}_{K}$.  By Lemmas~\ref{gens}
and~\ref{Wed}, it suffices to prove that $s-t$ is a nilpotent
element of $KS$ for all $s,t \in S$.

We first make some reductions. By Proposition~\ref{localgroups} that
there is an integer $n$ such that every product of at least $n$
elements in $S$ belongs to its minimal ideal $U$. In particular,
for all $s,t \in S$, $(s-t)^n$ belongs to $KU\cap \omega KS=\omega
KU$. Thus it suffices to show $\omega KU$ is nilpotent. In other
words, we may assume without loss of generality that $S$ is simple
with maximal subgroups in $\pv G_{K}$ and we change notation
accordingly.

Since $S$ is a simple semigroup, it constitutes a single $\J$-class
whence, as observed in Section~2, for every two $s,t \in S$, there
exists an idempotent $e\in S$ such that $s\R e$ and $e \L t$. Then
$$s-t = (s-e) + (e-t).$$ So the augmentation ideal of $S$ is
generated as a vector space by differences of elements in either
the same $\mathcal{R}$-class or the same $\mathcal{L}$-class,
with one of them being an idempotent.

Assume that $s\R e$ or $s\L e$ and let $f=f^2$ be the idempotent in
the $\mathcal{H}$-class of $s$. Then $$s-e =(s-f)+(f-e).$$ Thus we see
that the augmentation ideal is generated as a vector space by elements
that are either the difference of an element and the idempotent in its
$\mathcal{H}$-class or the difference of two idempotents in the
same $\mathcal{R}$-class or the same $\mathcal{L}$-class.

Consider an element $s-f$ where $f^2=f\H s$. If $\mathrm{char}\,K=0$,
then $s=f$, since $\pv G_{K}$ is the trivial variety and there is nothing
to prove. If $\mathrm{char}\,K = p$, then there is an $n$ such that
$s^{p^{n}}=f$ and since $s$ commutes with $f$, we have
$$(s-f)^{p^{n}}= s^{p^{n}}-f=0,$$ so in all cases $s-f$ is a
nilpotent element.

Now consider an element $f-e$ where $e$ and $f$ are idempotents and
either $f\L e$ or $f\R e$. Then
$$(f-e)^2=f^2-ef-fe+e^2=0$$ by \eqref{leftrightidentity}.

Therefore, the augmentation ideal is generated as a vector space
by nilpotent elements and we have shown $\omega KS$ is nilpotent.

The converse is a consequence of Lemma~\ref{toalgebra} with
$A=KS$, $B=K$ and $\p$ the augmentation map.
\end{proof}

\begin{Thm}\label{LG=nilpotent}
Let $\p:S \rightarrow T$ be a morphism of finite semigroups. Then
$\p$ is an $\pv{LG}_{K}$-morphism if and only if $\ker
\overline{\p}$ is a nilpotent ideal of $KS$.
\end{Thm}

\begin{proof}
Sufficiency is immediate from Lemma~\ref{toalgebra}.  For
necessity, suppose $\p$ is an $\pv{LG}_{K}$-morphism.  Then by
Lemmas~\ref{gens} and~\ref{Wed}, it suffices to prove that $s_1
-s_2$ is a nilpotent element of $KS$ for each $s_1,s_2 \in S$ with
$s_{1}\p = s_2\p$.

 Let $n$ be an integer such that $(s_1\p)^n$ is an idempotent
$f$ of $T$. Since $\p$ is an $\pv{LG}_{K}$-morphism, $U=f\pinv$ is
in $\pv {LG}_{K}$.  Also any product involving $n$ elements of the
set $\{s_1,s_2\}$ belongs to $U$.  Therefore, $$(s_1 - s_2)^n\in
KU\cap \omega KS = \omega KU$$ and is hence nilpotent by
Proposition~\ref{simp}.  It follows that $s_1-s_2$ is nilpotent,
as desired.
\end{proof}

Theorem \ref{LG=nilpotent} is a semigroup theorist's version of a
classical and central result of the theory of finite dimensional
algebras and holds in an appropriate sense for all finite
dimensional algebras. Indeed, it has been known since the early
1900's that if $A$ is a finite dimensional algebra and $N$ is a
nilpotent ideal of $A$, then every idempotent of $A/N$ lifts to an
idempotent of $A$. Furthermore, if we assume that algebras have
identity elements then two lifts of an idempotent in $A/N$ are
conjugate by an element of the group of units of $A$ of the form
$g=1+n, n \in N$ and more generally, a conjugacy class of
idempotents of $A/N$ can be lifted to a single conjugacy class of
idempotents in $A$. Lastly, two idempotents $e,f$ of a finite
dimensional algebra $A$ are conjugate by an element of the group
of units of $A$ if and only if $e \J f$ in the multiplicative
monoid of $A$. Putting this all together, it can be shown that
considered as a morphism between multiplicative monoids, the
morphism from $A$ to $A/N$ is an $\pv {LG}$-morphism.

\begin{Thm}\label{rhodesrad1}
The Rhodes radical of a finite semigroup $S$ over a field $K$ is
the largest $\pv{LG}_{K}$-congruence on $S$.
\end{Thm}

\begin{proof}
Since the map $KS\to KS/\mathrm{Rad}(KS)$ has nilpotent kernel,
Lemma~\ref{toalgebra} shows that $\Radf K S$ is an
$\pv{LG}_{K}$-congruence.  If $\p:S\onto T$ is any
$\pv{LG}_{K}$-morphism, then $\ov {\p}:KS\onto KT$ has nilpotent
kernel by Theorem~\ref{LG=nilpotent}, whence $\ker \ov \p\subseteq
\mathrm {Rad}(KS)$. Thus if $s_1\p = s_2\p$, then $s_1-s_2\in \ker \ov
\p\subseteq \mathrm{Rad}(KS)$, showing that $(s_1,s_2)\in \Radf K S$,
as desired.
\end{proof}

As a consequence we now give a simpler proof of some results of
Krohn--Rhodes--Tilson~\cite{Arbib,folleyR,folleyT}.

\begin{Lemma}\label{Radonto}
Let $\p:S\onto T$ be a surjective morphism of finite semigroups.
Then $\p$ induces a surjective morphism $\til \p:S/\Radf K S\onto
T/\Radf K T$.
\end{Lemma}
\begin{proof}
Clearly $\mathrm{Rad}(KS)\ov {\p}$ is a nilpotent ideal of $KT$
and hence contained in $\mathrm{Rad}(KT)$. Therefore a morphism
$KS/\mathrm{Rad}(KS)\onto KT/\mathrm{Rad}(KT)$ is well defined.
\end{proof}

\begin{Thm}\label{malcevbymorph}
Let \pv V be a variety of finite semigroups and $S$ a finite
semigroup. Then the following are equivalent:
\begin{enumerate}
\item $S\in\pv {LG}_{K}\malce \pv V$;
\item $S/\Radf K S\in \pv V$;
\item There is an $\pv {LG}_{K}$-morphism $\p:S\to T$ with
$T\in \pv V$.
\end{enumerate}
\end{Thm}

\begin{proof}
Since $\Radf K S$ is a $\pv {LG}_{K}$-congruence, (2) implies
(3).  Clearly (3) implies (1).  For (1) implies (2), suppose
$\p:T\onto S$ and $\psi:T\onto U$ are surjective morphisms with
$U\in \pv V$ and $\psi$ is an $\pv {LG}_{K}$-morphism. Then, by
Theorem~\ref{rhodesrad1}, the canonical morphism $T\onto T/\Radf K T$
factors through $\psi$ and so  $T/\Radf K T$ is a quotient of $U$ and
hence belongs to \pv V. It now follows from Lemma~\ref{Radonto} that
$S/\Radf K S\in \pv V$.
\end{proof}

The central question about a variety of finite semigroups is
usually the decidability of its membership problem. We say that a
variety \pv V is said to have \emph{decidable membership} if there
exists an algorithm to recognize whether a given finite semigroup
$S$ belongs to \pv V. The above results imply that the Mal'cev
product $\pv {LG}_{K}\malce \pv V$ has decidable membership
whenever the variety \pv V has. Indeed, given a finite semigroup
$S$, one effectively constructs its Rhodes radical $\Radf K S$ as
the largest $\pv{LG}_{K}$-congruence on $S$ and then verifies,
using decidability of membership in \pv V, the condition (2) of
Theorem~\ref{malcevbymorph}. This observation is important because
Mal'cev products of decidable pseudovarieties need not be
decidable in general. See~\cite{Rhodesundec,KBundec}.

\subsection{Rhodes Radical via Semigroup Theory}

We now indicate how to prove Theorem~\ref{rhodesrad1} using
semigroup representation theory.  Here we use the characterization
of the Rhodes radical as the intersection of the congruences
corresponding to all irreducible representations of $S$ over $K$.
This method will give us an explicit description of $\Radf K S$.

Krohn and Rhodes introduced the notion of a generalized group
mapping semigroup in~\cite{KR2}.  A semigroup $S$ is called
\emph{generalized group mapping}~\cite{Arbib,KR2} (GGM) if it has
a ($0$-)minimal ideal $I$ on which it acts faithfully on both the
left and right by left and right multiplication respectively. This
ideal $I$ is uniquely determined and is of the form
$I=J(\cup\{0\})$ where $J$ is a regular $\J$-class. We shall call
$I$ the \emph{apex} of $S$, written $\apex S$. We aim to show that
finite irreducible matrix semigroups are generalized group
mapping.

The following result was stated by Rhodes for the case of the field
of complex numbers~\cite{rhodeschar,RhodesZalc} but holds true in general.
Our proof for the general case uses the results of Munn and
Ponizovsky~\cite{CP,RhodesZalc}.

\begin{Thm}\label{thm:irredGGM}
Let $K$ be a field, $V$ be a finite dimensional $K$-vector space and
$S\leq \End K V$ be a finite, irreducible subsemigroup.  Then $S$ is
generalized group mapping.
\end{Thm}
\begin{proof}
If $S$ is the trivial semigroup, then it is clearly generalized
group mapping.  So we may assume $S$ is non-trivial.   Let $I$ be
a $0$-minimal ideal of $S$; if $S$ has no zero, take $I$ to be the
minimal ideal.  It is shown in~\cite[Theorem 5.33]{CP} that the
identity of $\End K V$ is a linear combination of elements of $I$.
We shall provide a proof of this for the sake of completeness.  It
will then immediately follow that $S$ acts faithfully on both the
left and right of $I$ by left and right multiplication.

The proof proceeds in several steps.  Let $I^{\perp} = \{v\in V\mid
vI=0\}$.  We first show that $I^{\perp} =0$.  To do this, we begin
by showing that $I^{\perp}$ is $S$-invariant.  Indeed, if $s\in S$,
$t\in I$ and $v\in I^{\perp}$, then using that $st\in I$, we
have$$(vs)t = v(st)=0,$$ showing that $vs\in I^{\perp}$.  Since
$I\neq \{0\}$, we cannot have $I^{\perp} = V$; thus $I^{\perp}=0$ by
irreducibility of $S$.

Next we show that $I$ itself is irreducible.  Let $\{0\}\neq W\leq
V$ be an $I$-invariant subspace.  Let $W_0 = \textrm{Span}\{wt\mid
w\in W,\ t\in I\}$. Notice that $W_0\subseteq W$.  If $w\in W$,
$t\in I$ and $s\in S$, then $(wt)s=w(ts)\in W$ since $ts\in I$ and
$W$ is $I$-invariant. Hence $W_0$ is $S$-invariant and so $W_0$ is
either $\{0\}$ or $V$. Since $I^{\perp}=0$, we cannot have that
$W_0=\{0\}$ and so $W\supseteq W_0=V$ establishing that $W=V$.  We
conclude that $I$ is irreducible.

Let $A$ be the $K$-span of $I$ inside of $\End K V$.  Then $A$ is an
irreducible algebra acting on $V$ and hence is a simple algebra by a
well-known result of Burnside \cite{Lam,CP}.  Thus $A$ has an
identity element $e$ by Wedderburn's theorem.  But $e$ commutes with
the irreducible semigroup $I$ and hence, by Schur's lemma, is
non-singular.  But the only non-singular idempotent endomorphism of
$V$ is the identity map and so the identity map belongs to $A$, the
linear span of $I$.
\end{proof}

\begin{Cor}
A finite irreducible subsemigroup of $M_n(K)$ has a unique 0-minimal
ideal, which is regular.
\end{Cor}

We recall some notions and results of Krohn and Rhodes.  The reader
is referred to~\cite{Arbib} for details.  Fix a finite semigroup $S$.
Choose for each regular $\J$-class $J$ a fixed maximal subgroup $G_J$.

\begin{Prop}[{\cite[Fact 7.2.1]{Arbib}}]\label{Prop:Jlifts}
Let $\p:S\onto T$ be a surjective morphism.  Let $J'$ be a $\J$-class
of $T$ and let $J$ be a $\j$-minimal $\J$-class of $S$ with $J\p\cap
J'\neq \emptyset$.  Then $J\p =J'$.  Moreover if $J'$ is regular,
then $J$ is unique and regular, and the images of the maximal
subgroups of $J$ are precisely the maximal subgroups of $J'$.
\end{Prop}

If $T$ is GGM and $J'=\apex T\setminus 0$, then we shall call the
$\J$-class $J$ of the above proposition the \emph{apex} of $\p$,
denoted $\apex {\p}$.  Let $K_{\p}$ be the group theoretic kernel
of $\p|_{G_{\apex {\p}}}$. We call $K_{\p}$ the \emph{kernel} of
$\p$. Krohn and Rhodes showed~\cite{Arbib} that $\p$ is completely
determined by its apex and kernel.

Let $J$ be a regular $\J$-class of $S$ and $N\lhd G_J$ be a normal
subgroup. We denote by $R_a$, $a\in A$, the $\R$-classes of $J$
and by $L_b$, $b\in B$, the $\L$-classes of~$J$.  Suppose that
$G_J=R_1\cap L_1$. For each $a\in A, b\in B$, choose according to
Green's Lemma~\cite{CP} $r_a\in J$ such that $s\mapsto r_as$ is a
bijection $R_a\to R_1$ and $l_b\in J$ such that $s\mapsto sl_b$ is
a bijection $L_b\to L_1$.  With this notation if $H_{ab}=R_a\cap
L_b$, then $s\mapsto r_asl_b$ is a bijection $H_{ab}\to G_J$.

We define a congruence by $s\equiv_{(J,G_J,N)}t$ if and only if,
for all $x,y\in J$,
\begin{equation}\label{ggmdef1}xsy\in J\iff xty\in J
\end{equation}
and, in the case where $xsy\in J$, if $x\in R_a$ and $y\in L_b$, then
\begin{equation}\label{ggmdef2}
r_axsyl_bN=r_axtyl_bN.
\end{equation}
  The quotient $S/\!{\equiv_{(J,G_J,N)}}$ is
denoted $\GGM (J,G_J,N)$~\cite{Arbib}.  The following result is
the content of~\cite[Proposition 8.3.28, Remark 8.3.29]{Arbib}.

\begin{Thm}\label{Thm:ggmthm}
Let $S$ be a finite semigroup and $J$ a regular $\J$-class with
maximal subgroup $G_J$. Suppose $\p:S\onto T$ is a surjective
morphism with $T$ a generalized group mapping semigroup.  Let
$J=\apex {\p}$ and let  $K_{\p}$ be the kernel.  Then the
congruence associated to  $\p$ is $\equiv_{(J,G_J,K_{\p})}$.  In
particular, $T\cong
  \GGM (J,G_J,K_{\p})$.
\end{Thm}

It follows from the above theorem that the definition of $\GGM
(J,G_J,K_{\p})$ doesn't depend on the choices made. The following
result from~\cite{Arbib} is an immediate consequence of the
definition of $\equiv_{(J,G_J,N)}$.

\begin{Cor}\label{factor}
Let $\p_1:S\onto T_1$, $\p_2:S\onto T_2$ be surjective morphisms
to generalized group mapping semigroups with $\apex
{\p_1}=\apex{\p_2}$. Then $\p_2$ factors through $\p_1$ if and
only if $K_{\p_1}\leq K_{\p_2}$.
\end{Cor}

We shall need the following fundamental result on semigroup
representations, due to  Clifford, Munn and Ponizovsky, which is one
of the main results of~\cite[Chapter 5]{CP} (see
also~\cite{RhodesZalc}).

\begin{Thm}\label{Thm:reptheory}
Let $S$ be a finite semigroup, $J$ a regular $\J$-class of $S$ and
$G_J$ a maximal subgroup of $J$.  Then any irreducible
representation \mbox{$\rho:G_J\to GL(V)$} can be extended uniquely
to an irreducible representation of $S$ with apex $J$.
\end{Thm}

It is proved in~\cite{CP,RhodesZalc} that every irreducible
representation of a finite semigroup $S$ is obtained by extending an
irreducible representation of a maximal subgroup $G_J$ for some
regular $\J$-class $J$, although we shall not need this result.

We are now ready to prove Theorem~\ref{rhodesrad1} via
representation theoretic means. First we need the following
classical result, which is a consequence of Maschke's theorem and
Clifford's theorem from finite group representation theory,
handling the group case~\cite[Corollary 8.6]{Lam}\footnote{We
thank John Dixon for pointing this result out to us.}. If $G$ is a
finite group, we define $G_{K}$, called the \emph{unipotent
radical} of $G$, to be the largest normal subgroup of $G$
belonging to $\pv {G}_{K}$.

\begin{Thm}\label{Thm:groupcase}
Let $G$ be a group and $K$ be a field.  Then $\Radf K G$ is the
congruence whose classes are the cosets of $G_{K}$.
\end{Thm}

We remark that Theorem~\ref{Thm:groupcase} also follows from our
first proof of Theorem~\ref{rhodesrad1} since the largest $\bf
LG_K$ congruence on a finite group clearly has kernel ${G}_{K}$.

\begin{Thm}\label{Thm:Rhodesradical}
Let $S$ be a semigroup and $K$ be a field.  Then $\Radf K S$ is
the congruence associated to the direct sum over all regular
$\J$-classes $J$ of the maps
\begin{equation}\label{radicalequation}
S\onto \GGM (J,G_J,(G_J)_{K}).
\end{equation}
\end{Thm}
\begin{proof}
Let $\sim$ be the congruence associated to the direct sum of the
maps \eqref{radicalequation}.  Let $\p$ be an irreducible
representation of $S$ with apex $J=\apex {\p}$. Then, by
Theorem~\ref{Thm:groupcase}, $(G_J)_{K}\leq K_{\p}$ and so, by
Corollary~\ref{factor}, $\p$ factors through $S\onto \GGM
(J,G_J,(G_J)_{K})$.  Thus $\sim\ \subseteq \Radf K S$.

For the reverse inclusion, it suffices to show that the congruence
associated to each map $S\onto \GGM (J,G_J,(G_J)_{K})$ can be
realized by a direct sum of irreducible representations.  Fix a
regular $\J$-class $J$ and let $\{N_i\}$ be a collection of normal
subgroups of $G_J$. Then the congruence associated to the direct
sum of the maps $S\onto \GGM (J,G_J,N_i)$ is the congruence
associated to $S\onto \GGM (J,G_J,\bigcap N_i)$. In particular,
consider the collection $\{K_{\p}\}$ where $\p$ is an irreducible
representation of $S$ with apex $J$.  Then, by
Theorem~\ref{Thm:reptheory}, the $K_{\p}$ run over all kernels of
irreducible representations of $G_J$, so by
Theorem~\ref{Thm:groupcase}, we obtain $\bigcap K_{\p} =
(G_J)_{K}$. The theorem now follows.
\end{proof}

Notice that this theorem allows for an explicit determination of
$\Radf K S$ via \eqref{ggmdef1} and \eqref{ggmdef2}.  The fact
that the above congruence is the largest $\pv
{LG}_{K}$-congruence is contained
in~\cite{Arbib,folleyR,folleyT}.

\section{Applications to Diagonalizability and Triangularizability}
Our first application of the Rhodes radical is to the question of
diagonalizablity and triangularizability of finite semigroups.
In~\cite{amv}, three of the authors characterized the varieties of
finite semigroups that can be (uni)triangularized over finite
fields. Using our techniques, we give a shorter, more conceptual
proof that works over a general field.

Let $K$ be a field.  Define $\pv {Ab}_K$ to be the variety of finite
Abelian groups generated by all finite subgroups of $K^*$.  It is well
known that any finite subgroup $G$ of $K^*$ is cyclic and is the set of
roots of $x^{|G|}-1$.  Moreover, there is a cyclic subgroup of $K^*$
of order $m$ if and only if $x^m-1$ splits into distinct linear
factors over $K$.  It is not hard to see that if $x^e-1$ and $x^f-1$
split into distinct linear factors, then so does
$x^{\mathrm{lcm}(e,f)}-1$.   Also if $x^e-1$ splits into linear
factors, then so does $x^d-1$ for any divisor $d$ of $e$.  Hence $\pv
{Ab}_K$ can be described as the variety of all
finite Abelian groups whose exponent $e$ has the property that $x^e-1$
splits into $e$ distinct linear factors over $K$.  We remark that if
the characteristic of $K$ is $p>0$, then $e$ and $p$ must be
relatively prime for this to happen.  If $K$ is algebraically closed
of characteristic $0$, then $\pv {Ab}_K=\pv {Ab}$.  If $K$ is
algebraically closed of characteristic $p>0$, then $\pv {Ab}_K$
consists of all finite Abelian $p'$-groups, that is of all finite
Abelian groups whose orders are relatively prime to $p$.

If $\pv H$ is a variety of finite groups, then the elements of the
variety $\pv H\malce \pv {Sl}$ are referred to as
\emph{semilattices of groups from \pv H}.  Such semigroups are
naturally ``graded'' by a semilattice in such a way that the
homogeneous components (which are the $\mathcal H$-classes) are
groups from $\pv H$.  It turns out that $\pv H\malce \pv {Sl}$ is
the varietal join $\pv {Sl}\vee \pv H$. See \cite{CP,almbook} for
more details. The following exercise in Linear Algebra captures
diagonalizability.

\begin{Thm}\label{Thm:diag}
Let $K$ be a field and $S$ a finite semigroup.  Then the following are
equivalent:
\begin{enumerate}
\item $S$ is commutative and satisfies an identity $x^{m+1}=x$ where
  $x^m-1$ splits into distinct linear factors over $K$;

\item $S$ is a semilattice of Abelian groups from $\pv {Ab}_K$;

\item Every representation of $S$ is diagonalizable;

%\item $KS\cong K^n$ for $n=|S^{1}|$;

\item $S$ is isomorphic to a subsemigroup of $K^n$ for
$n=|S^{1}|$.
\end{enumerate}
\end{Thm}
\begin{proof}
The equivalence of (1) and (2) follows from Clifford's
Theorem~\cite[Theorem 4.11]{CP}.  For
 (1) implies (3), suppose $\rho:S\to \End K V$ is a representation.
 Since $S$ satisfies $x^{m+1}=x$, we must have that $s\rho$ satisfies
 $x(x^m-1)=0$.  It
 follows that the minimal polynomial of $s\rho$ for any $s\in S$ has
 distinct roots and
 splits
 over $K$. Hence $s\rho$ is diagonalizable for all $s\in S$. To show
 that $\rho$ is diagonalizable we induct on the degree of the
 representation.  If $\rho$ is of degree one, then clearly it is
 diagonalizable.  If $S\rho$ is contained in the scalar matrices, then
 we may also deduce that the representation is diagonalizable.
 Otherwise, there is an element $s\in S$ such that $s\rho$ is not a
 scalar matrix.  Since $s\rho$ is diagonalizable, we can write
 $V=\bigoplus_{\lambda\in \mathrm{Spec}(s\rho)} E_{\lambda}$ where
 $E_{\lambda}$ is the eigenspace of $\lambda$.  We claim that
 $E_{\lambda}$ is $S$-invariant. Indeed, if $t\in S$ and $u\in
 E_{\lambda}$, then $$ut\rho s\rho = us\rho t\rho = \lambda ut\rho$$
 and so $ut\rho\in E_{\lambda}$.  Since $s\rho$ is not a scalar, each
 $E_{\lambda}$ has smaller dimension and so the restriction of $\rho$ can be diagonalized by
 induction.  Thus we have diagonalized $\rho$.

(3) implies (4) follows immediately from considering the right
regular representation of $S$, that is by having $S$ act on $S^1$
by right multiplication and extending linearly. To show that (4)
implies (1), first observe that $S$ embeds in a direct product of
finite subsemigroups of $K$.   A finite subsemigroup of $K$ is
commutative and must satisfy an identity of the form
$x^i=x^{p+i}$, for some $i\geq 0, p >0$. Since $K$ is a field, we
deduce that $S$
     satisfies an identity of the form $x^{m+1}=x$ for some $m\geq 0$.    If $m$
     is minimum, then $x^m-1$ splits into distinct linear factors, as discussed above.
 This completes the proof.
\end{proof}

With a little more work, we can improve condition (4) in the
previous theorem a bit. It is easy to prove that any semilattice
of groups $S$ acts faithfully by right multiplication on $S$
considered as a set, even if $S$ does not have an identity. Thus
we can replace condition (4) by $S$ is isomorphic to a
subsemigroup of $K^n$ for $n=|S|$. From this it immediately
follows that the finite semigroups described in Theorem 4.1 are
precisely the finite semigroups $S$ such that $KS \cong K^{n}$
where $n=|S|$. Details are left to the reader.

 The semigroups satisfying the
conditions of the above theorem for diagonalizability  form a
variety of finite semigroups that we denote by $\pv D_K$. It is
precisely the variety $\pv{Ab}_K\malce \pv {Sl}$ by Clifford's
Theorem~\cite[Theorem 4.11]{CP}.  It is the varietal join $\pv
{Sl}\vee \pv {Ab}_K$, which can be seen from the above theorem.
For example, if $K=\mathbb{C}$, then $\pv D_K$ consists precisely
of semilattices of Abelian groups (i.e.\ commutative inverse
semigroups).  If $K=\mathbb{F}_q$, the finite field of $q$
elements, then $\pv D_K$ consists of semilattices of Abelian
groups with exponent dividing $q-1$.  If $K$ is the algebraic
closure of $\mathbb{F}_p$, then $\pv D_K$ consists of semilattices
of Abelian $p'$-groups.

We shall call a matrix
\emph{unidiagonal} if it is diagonal and its entries are contained
in $\{0,1\}$.  We have the following analogue of Theorem
\ref{Thm:diag} whose proof we leave to the reader.

\begin{Thm}\label{Thm:unidiag}
Let $K$ be a field and $S$ a finite semigroup.  Then the following are
equivalent:
\begin{enumerate}
\item $S$ is a semilattice;
\item Every representation of $S$ is unidiagonalizable;
\item $KS \cong K^n$ for $n=|S|$ and the image
of $S$ under this isomorphism is contained in $\{0,1\}^n$;

\item $S$ is isomorphic to a subsemigroup of $\{0,1\}^n$ for some $n$.
\end{enumerate}
\end{Thm}

The above theorem shows that \pv {Sl} is the variety of finite
unidiagonalizable semigroups.

Let $K$ be a field.  Let $T_n(K)$ denote the semigroup of upper
triangular $n\times n$ matrices over $K$.  Recall that a finite
dimensional $K$-algebra $A$ is called \emph{basic} if
$A/\mathrm{Rad}(A)$ is commutative.  If $A/Rad(A)\cong K^n$ for
some $n$, then $A$ is called a \emph{split basic $K$-algebra}.

\begin{Thm}\label{Thm:triangularizability}
Let $S$ be a finite semigroup and $K$ a field. Then the following
are equivalent.
\begin{enumerate}
\item $S\in \pv {LG}_{K}\malce \pv D_K$ ;
\item $S/\Radf K S\in \pv D_K$ ;
\item $KS/\mathrm{Rad}(KS) \cong K^m$ for some $m$;
\item $KS$ is a split basic $K$-algebra;
\item Every irreducible representation of $S$ over $K$ has degree one;
\item Every representation of $S$ is triangularizable;
\item $S\leq T_n(K)$, $n=|S^1|$;
\item $S\leq T_m(K)$, some $m$.
\end{enumerate}
\end{Thm}
\begin{proof}
We have already seen the equivalence of (1) and (2).  For (2)
implies (3), let $T=S/\Radf K S$.  Then $KT\onto
KS/\mathrm{Rad}(KS)$.  By Theorem~\ref{Thm:diag} $KT=K^{|T|}$. Hence
$KS/\mathrm{Rad}(KS)$ is a direct product of copies of $K$.  The
equivalence of (3) and (4) is the definition.

The implication (3) implies (5) follows immediately from the
Wedderburn theory, since the only irreducible representations of a
direct product of fields are the projections. For (5) implies (6),
let $\p:S\to M_m(K)$ be a representation. Then by choosing a
composition series for the right $KS$-module associated to $\p$,
we can put $S\p$ in block upper triangular form where the diagonal
blocks are irreducible representations or the zero representation.
But since all such are of degree one, we conclude that $S\p$ has
been brought to triangular form.

One establishes (6) implies (7) by considering the regular
representation of $S$.  That (7) implies (8) is trivial. For (8)
implies (1), observe that the projection $\p$ from $T_m(K)$ to the
diagonal is an algebra homomorphism with nilpotent kernel. Thus
$\p|_S$ is an $\pv{LG}_{K}$-morphism by Lemma~\ref{toalgebra} and so
(1) follows from Theorem \ref{Thm:diag}.
\end{proof}

 Let $UT_n(K)$ denote the semigroup
of upper unitriangular $n\times n$ matrices over $K$, where by
unitriangular we mean triangular with only $1$'s and $0$'s on the
diagonal.
 By a \emph{trivial representation} of $S$, we mean a
homomorphism $\p:S\to \{0,1\}$.  The following theorem is proved
similarly to the above theorem.  We omit the proof.

\begin{Thm}\label{Thm:unitriangularizability}
Let $S$ be a finite semigroup and $K$ a field. Then the following
are equivalent.
\begin{enumerate}
\item $S\in \pv {LG}_{K}\malce \pv
{Sl}$;
\item $S/\Radf K S\in \pv
{Sl}$;
\item $KS/\mathrm{Rad}(KS) \cong K^m$ for some $m$ and the image of
  $S$ is contained in
$\{0,1\}^m$;
\item Every irreducible representation of $S$ over $K$ is
trivial;
\item Every representation of $S$ is unitriangularizable;
\item $S\leq UT_n(K)$, $n=|S^1|$;
\item $S\leq UT_m(K)$, some $m$.
\end{enumerate}
\end{Thm}

Notice that unitriangularizability depends only on the
characteristic and not the field.  The proofs of condition (6) in
the above theorems show that a [uni]triangularizable monoid can be
realized as a submonoid of $S\leq T_n(K)$ $[S\leq UT_n(K)]$ and a
[uni]triangular group can be realized as a subgroup of $T^*_n(K)$
[$UT_n^*(K)$] (where here $*$ denotes the group of units of a
monoid). We remark that if a finite semigroup $S$ is
triangularizable over the algebraic closure $\ov K$ of $K$, then
it is triangularizable over a finite extension of $K$. Indeed, $S$
can be faithfully represented in $T_n(K)$ ($UT_n(K)$), where
$n=|S^1|$. Since only finitely many entries appear amongst the
entries of $S$, we can just take the extension field generated by
these entries.  The same remarks apply to diagonalization.

We now determine the above varieties.   Recall that if \pv H is a
variety of finite groups, then $\ov {\pv H}$ denotes the variety
of finite semigroups all of whose subgroups belong to \pv H.
Usually $\ov {\pv I}$ is denoted \pv A (for aperiodic). If $\pv V$
is a variety of finite semigroups, then $\pv {DV}$ is the variety
of semigroups whose regular $\J$-classes are subsemigroups that
belong to $\pv V$.  If \pv H is a variety of finite groups, then
$\pv {Sl}\vee \pv H$ is the variety of semilattices of groups from
\pv H \cite{almbook}.  We denote by \pv O the variety of finite
orthodox simple semigroups.  A simple semigroup $S$ is
\emph{orthodox} if $E(S)$ is a subsemigroup. If \pv V is a variety
of finite semigroups, \pv {EV} is the variety of finite semigroups
$S$ such that $E(S)$ generates a subsemigroup in \pv V.

To handle the case of characteristic zero,  we need a result that
can easily be verified by direct calculations with generalized group
mapping congruences.  Since a syntactic proof can be found in
\cite[Corollary 3.3]{AlmeidaDA} we skip the proof.

\begin{Lemma}\label{DOlemma}
Let \pv H be a variety of finite groups, then $$\pv {LI}\malce (\pv
{Sl}\vee \pv H) = \pv {DO}\cap \ov{\pv H}.$$
\end{Lemma}

In particular we obtain the following corollary:

\begin{Cor}\label{DAisforever}
The variety of unitriangularizable semigroups in characteristic
zero is \pv {DA}. The variety of triangularizable semigroups over
a field $K$ of characteristic zero is $\pv {DO}\cap
\ov{\pv{Ab}_K}$. In particular, the variety of triangularizable
semigroups over an algebraically closed field of characteristic
zero is $\pv {DO}\cap \ov{\pv {Ab}}$.
\end{Cor}

Of course \pv{DA} and $\pv {DO}\cap \ov {\pv {Ab}}$ are decidable
varieties.  In general, decidability of $\pv {DO}\cap
\ov{\pv{Ab}_K}$ depends on $K$.  Notice that \pv {DA} contains all
finite bands, that is, all finite idempotent semigroups. The
triangularizability of bands can be found in the work of
\cite{Brown1,Brown2}. Corollary \ref{DAisforever} is useful for
computing spectra of random walks on semigroups in \pv{DA} or $\pv
{DO}\cap \ov{\pv{Ab}_K}$\cite{Mobius}. In particular, some famous
Markov chains, such as the Tsetlin library, arise as random walks
on bands \cite{BHR,Brown1,Brown2}. Another consequence of
Corollary \ref{DAisforever} is that the semigroup algebra of a
finite semigroup $S$ is split basic over the reals if and only if
$S\in \pv{DO}$ and every subgroup of $S$ has exponent two.

We now turn to the case of characteristic $p$.

\begin{Lemma}\label{pversionofDOlemma}
Let $p$ be a prime and let \pv H be a variety of finite
$p'$-groups.  Then
\begin{equation}\label{eq:DOlem}
\pv {LG}_p\malce (\pv {Sl}\vee \pv H) = \pv D(\ov{\pv {G}_p\malce
\pv H})\cap \pv E\ov{\pv {G}_p}.
\end{equation}
\end{Lemma}
\begin{proof}
To see that the left hand side of \eqref{eq:DOlem} is contained in
the right hand side, suppose $S\in \pv{LG}_p\malce (\pv {Sl}\vee \pv
H)$. Let $T = S/\mathrm{Rad}_{\FF_p}(S)$ and let $\p:S\to T$ be the
canonical homomorphism. Then $\p$ is an $\pv{LG}_p$-morphism and
$T\in \pv {Sl}\vee \pv H$ by Theorem~\ref{malcevbymorph}.  Hence if
$J$ is a regular $\J$-class of $T$, then $J\pinv$ is a nilpotent
extension of a simple semigroup by Proposition~\ref{localgroups}. It
easily follows that $S\in
\pv{DS}$ (since regular $\J$-classes are mapped into regular
$\J$-classes). Suppose $G$ is a subgroup of $S$. Then $G\p\in \pv H$
and $\ker \p|_G\in \pv {G}_p$ since $\p$ is an $\pv
{LG}_p$-morphism. We conclude $S\in \pv D(\ov{\pv {G}_p\malce \pv
H})$.  Let $J$ be a regular $\J$-class of $S$. Let $E(J)$ be the
idempotents of $J$. Then $E(J)\p$ is the unique idempotent $f$ of
the $\J$-class $J\p$ of $T$ (since $T$ is a semilattice of groups).
Hence $\langle E(J)\rangle \p = f$. Since $f\pinv \in \pv {LG}_p$ it
follows that every maximal subgroup of $\langle E(J)\rangle$ belongs
to $\pv {G}_p$. This shows that $S\in \pv E\ov {\pv {G}_p}$.  This
establishes the inclusion from left to right in \eqref{eq:DOlem}.

For the reverse inclusion, it suffices to show that if $J$ is a
regular $\J$-class of a finite semigroup $S$ in the right hand side
of \eqref{eq:DOlem}, then $$\GGM (J,G_J, (G_J)_{\FF_p})\in \pv
{Sl}\vee \pv H.$$ First note that since $\pv H$ consists of
$p'$-groups, $G_J\in \pv G_p\malce \pv H$ means precisely that $G_J$
has a normal $p$-Sylow subgroup $N$ and that $G_J/N\in \pv H$.  We
remark that $N$ is the $p$-radical $(G_J)_{\FF_p}$. By the results
of~\cite{graham}, there is a Rees matrix representation
$\MM^0(G_J,A,B,C)$ of $J^0$ with the entries of $C$ generating the
maximal subgroup $K$ of the idempotent-generated subsemigroup. Since
$S\in \pv E\ov {\pv G_p}$, $K$ is a $p$-subgroup of $G_J$ and hence
contained in $N$. According to~\cite[8.2.22 Fact (e)]{Arbib}) to
obtain the image of $J$ in $\GGM (J,G_J,(G_J)_{\FF_p})$, we
project to $\MM^0(G_J/N,A,B,\ov C)$, where $\ov C$ is obtained from
$C$ by first reducing modulo $N$, and then identifying proportional
rows and columns.  But since the entries of $C$ belong to $N$, this
results in identifying all rows and columns and so the image of $J$
in $\GGM (J,G_J,(G_J)_{\FF_p})$ is simply $G_J/N$. Since $\GGM
(J,G_J,(G_J)_{\FF_p})$ acts faithfully on the right of its apex by
partially defined right translations and the only non-zero,
partially defined right translations of a group are zero and right
translations by elements of the group, we see that
$$\GGM (J,G_J,(G_J)_{\FF_p}) = G_J/N\ \text{or}\ (G_J/N)\cup
0$$ (depending on whether $J$ is the minimal ideal, or not).  Thus
$$\GGM (J,G_J,(G_J)_{\FF_p})\in \pv {Sl}\vee \pv H$$ as desired.
\end{proof}

Observing that extensions of $p$-groups by Abelian groups are the
same thing as extensions of $p$-groups by Abelian $p'$-groups, we
have the following corollary.

\begin{Cor}\label{cor:globaltriangularizability}
The variety of unitriangularizable semigroups over any field of
characteristic $p$ is $\pv D\ov {\pv G}_p$.  The variety of
triangularizable semigroups over a field $K$ of characteristic $p$
is $\pv D(\ov{\pv {G}_p\malce \pv {Ab}_K})\cap \pv E\ov{\pv G_p}$.
In particular, the variety of semigroups triangularizable over an
algebraically closed field of characteristic $p$ is $\pv D(\ov{\pv
{G}_p\malce \pv {Ab}})\cap \pv E\ov {\pv G_p}.$
\end{Cor}

In particular, commutative semigroups are triangularizable over
any characteristic. More precisely, every finite commutative
semigroup is triangularizable over some field of characteristic 0
and for some field of characteristic $p$ for each prime $p$. In
fact, the semigroups triangularizable over any characteristic are
precisely those in $\pv {DO}\cap \ov {\pv {Ab}}$. Pseudoidentities
for many of these varieties can be found in~\cite{amv}.  A method
of constructing pseudoidentities for $\pv {LI}\malce \pv V$ from
those of $\pv V$ can be found in~\cite{ReillyZ} and for $\pv
{LG}_p\malce \pv V$ from those of $\pv V$ can be found
in~\cite{newus}.

We now turn to characterize those finite semigroups whose
semigroup algebras are basic over a field $K$.  The case of split
basic $K$-algebras has already been handled in Theorem
\ref{Thm:triangularizability}.  Recall that a finite dimensional
$K$-algebra $A$ is called \emph{basic} if $A/\mathrm{Rad}(A)$ is
commutative, or equivalently, a direct product of fields. Since
$KS/\mathrm{Rad}(KS)$ is generated as an algebra by $S/\Radf K S$,
to be basic $S/\Radf K S$ must be a semilattice of Abelian groups
(embedding in a direct product of fields). Conversely, if $S/\Radf
K S$ is a semilattice of Abelian groups, then
$KS/\mathrm{Rad}(KS)$ (being generated by $S/\Radf K S$) must be a
commutative algebra. Thus we have proved:

\begin{Cor}
Let $S$ be a finite semigroup and $K$ a field.  Then $KS$ is a
basic algebra if and only if $$S\in \pv {LG}_{K}\malce (\pv
{Sl}\vee \pv {Ab}) = \begin{cases} \pv {DO}\cap \ov {\pv {Ab}} &
char\ K =0 \\ \pv D\ov{(\pv {G}_p\malce \pv {Ab})}\cap \pv
{E\ov{G}_p} & char\ K=p.\end{cases}$$
\end{Cor}

\section{Applications to Semigroup Decomposition Theory}
Our  next application of the Rhodes radical is to recover some
deep algebraic decomposition results of Rhodes--Tilson--Weil. For
the remainder of the paper we will deal with monoids and varieties
of finite monoids.

First of all we recall the definition of the two-sided semidirect
product of two monoids.  Let $M$ and $N$ be monoids and suppose
that $N$ has a bi-action on $M$ (that is commuting left and right
actions on $N$).  For convenience we write $M$ additively and $N$
multiplicatively although we assume no commutativity. Then the
\emph{two-sided semidirect product} $M\bowtie N$ consists of all
$2\times 2$ upper triangular matrices
$\begin{pmatrix} n & m\\
                  0 & n\end{pmatrix}$
with the usual matrix multiplication.  There is an obvious
projection to $N$ via the diagonal.  The variety generated
by two-sided semidirect products $M\bowtie N$ with $M\in \pv V$
and $N\in \pv W$ is denoted $\pv V\dsd \pv W$.

Rhodes and Tilson introduced in~\cite{Kernel} the kernel category
as a way to determine membership in $\pv V\dsd \pv W$. We restrict
ourselves to considering the kernel category of a morphism and to
a special case of the results of~\cite{Kernel} to avoid getting
technical. Let $\p:M\to N$ be a homomorphism.
Following~\cite{Kernel}, we define a category $\mathsf K_{\p}$,
called the \emph{kernel category} of $\p$.  The object set is
$N\times N$. The arrows are equivalence classes of triples
$(n_L,m,n_R)\in N\times M\times N$ where $(n_L,m,n_R):(n_L,m\p
n_R)\to (n_Lm\p,n_R)$ and two coterminal triples $(n_L,m,n_R)$ and
$(n_L,m',n_R)$ are identified if and only if
\mbox{$m_Lmm_R=m_Lm'm_R$} for all $m_L\in n_L\pinv, m_R\in
n_R\pinv$.  Composition is given by
$$[(n_L,m,m'\p n_R)][(n_Lm\p,m',n_R)] = [(n_L,mm',n_R)];$$ the
identity at $(n_L,n_R)$ is $[(n_L,1,n_R)]$.

We consider categories as partial algebras whose elements consist
of all of its arrows. If $C$ is a category and $c$ is an object of
$C$, then the collection of all arrows $C(c,c)$ from $c$ to itself
is a monoid called the {\emph {local monoid at $c$}}. It is clear
that if we add a new zero element to $C$, we obtain a semigroup
$C^{0}$ called the {\emph consolidation} of $C$. In $C^{0}$, the
identity element $e_c$ at $c$ is an idempotent and then it is easy
to see that $C(c,c)$ is isomorphic to the local monoid (in the
sense of our previous usage of that term in semigroup theory)
$e_{c}C(c,c)e_{c}$.

Let \pv V be a variety of finite monoids. A category $C$ is said
to be \emph{locally in $\pv V$} if each of the local monoids
$C(c,c)$ belongs to $\pv V$; we use the notation $C(c)$ as a
shorthand for $C(c,c)$. The collection of categories locally in
\pv V is denoted $\ell \pv V$. The following is an amalgamation of
results of~\cite{TilsonCat} and a special case of the results of
\cite{Kernel}.

\begin{Thm}
Let $M$ be a finite monoid, \pv H a non-trivial variety of finite
groups and \pv V a variety of finite monoids.  Then $M\in \pv H\dsd
\pv V$ if and only if there is a finite monoid $N$ mapping onto
$M$ that admits a morphism $\p:N\to V\in \pv V$ such that $\mathsf
K_{\p}\in \ell\pv H$.
\end{Thm}

Notice that $\ell \pv I$ is a variety of finite
categories~\cite{TilsonCat}, from which it easily follows that the
collection of finite monoids $M$ that are quotients of finite
monoids $N$ admitting a morphism to $\p:N\to V\in \pv V$ with
$\mathsf K_{\p}\in \ell \pv I$ is a variety of finite monoids,
which we denote $\ell \pv I\dsd \pv V$. This variety plays
an important role in language theory~\cite{Pinetal}, as we shall
see below.

Let $\pv {Ab}(p)$ denote the variety of finite Abelian groups of
exponent $p$, where $p$ is a prime. Our goal is to prove the following
two important cases of the results of Rhodes--Tilson--Weil~\cite{Kernel,MPS2}
(see also~\cite{Pinetal}).

\begin{Thm}\label{Thm:dbsd}
Let $\pv V$ be a variety of finite monoids and $p$ a prime.  Then
the smallest variety of finite monoids containing \pv V and closed
under the operations $\pv W\mapsto \ell \pv I\dsd \pv W$,
respectively $\pv W\mapsto \pv {Ab}(p)\dsd \pv W$, is $\pv
{LI}\malce \pv V$, respectively $\pv {LG}_p\malce \pv V$.
\end{Thm}

The original proof of Theorem~\ref{Thm:dbsd} is a case-by-case
analysis using Rhodes's classification of maximal proper
surjective morphisms~\cite{MPS,MPS2,Arbib}.  We give a conceptual
proof via representation theory.   First we make some preliminary
observations.

It is well known~\cite{Kernel,Arbib,folleyR,folleyT} that a
morphism is an $\pv{LI}$-morphism (respectively $\pv
{LG}_p$-morphism) if and only if it is injective on two element
semilattices and on subgroups (respectively on $p'$-subgroups). It
follows immediately that \pv {LI}-morphisms (respectively $\pv
{LG}_p$-morphisms) are closed under composition.  Thus if $K$ is a
field,
\begin{equation}\label{malcevclosed}
\pv {LG}_{K}\malce (\pv {LG}_{K}\malce \pv V)= \pv
{LG}_{K}\malce \pv V
\end{equation}

The following is well known (cf.~\cite{Kernel,MPS2,Pinetal}),
but we include the proof
for completeness.

\begin{Prop}\label{lVimpliesLV}
Let $\pv V$ be a variety of finite monoids and \mbox{$\p:M\to N$}
be a morphism with $\mathsf K_{\p}$ locally in $\pv V$.  Then $\p$
is an \pv{LV}-morphism.
\end{Prop}
\begin{proof}
Let $f\in E(N)$.  Set $M_f = f\pinv$ and let $m\in f\pinv$.  Then
$[(f,m,f)]:(f,f)\to (f,f)$ is an arrow of $\mathsf K_{\p}$.  Let
$e\in E(f\pinv)$ and define a map \mbox{$\psi:eM_fe\to
K_{\p}((f,f),(f,f))$} by $m\mapsto [(f,m,f)]$. Clearly this is a
morphism; we show it is injective. Suppose $m\psi =m'\psi$. Then
since $e\in f\pinv$, this implies $m =eme=em'e=m'$. Thus $eM_fe\in
\pv V$ and so $M_f\in \pv {LV}$, establishing that $\p$ is a \pv
{LV}-morphism.
\end{proof}

 Let $M_{m,r}(K)$ denote the collection of
$m\times r$ matrices over a field $K$. The following lemma will
afford us the decompositions needed for our proof of
Theorem~\ref{Thm:dbsd}.

\begin{Lemma}\label{lemma:decomp}
Let $K$ be a ring and $M\leq M_n(K)$ be a finite monoid of
block upper triangular matrices of the form $$\left\{\begin{pmatrix} A & B \\
                                                 0 &
C\end{pmatrix}\mid A\in M_m(K),\ B\in M_{m,r}(K),\ C\in
M_r(K)\right\}.$$  Let $N$ be the quotient of $M$ obtained by
projecting to the block diagonal and let $\p$ be the projection. Then
each local monoid of\/ $\mathsf K_{\p}$ embeds in the additive group
of $M_{m,r}(K)$.  In particular, if $K$ is a field, then
$$\mathsf K_{\p}\in
\begin{cases} \ell \pv I & char\ K=0 \\ \ell \pv {Ab}(p) & char\
K=p.\end{cases}$$
\end{Lemma}
\begin{proof}
Elements of $N$ are certain pairs $(A,C)$ with $A\in M_m(K)$ and
$C\in M_r(K)$.  Let $S=\mathsf K_{\p}(((X,Y),(U,V)))$. We define a
map $\psi:S\to M_{m,r}(K)$ as follows.   Given an arrow
\mbox{$a=[((X,Y),m,(U,V))]\in S$} with
$m=\begin{pmatrix} A & B \\
                                                 0 &
C\end{pmatrix}$, define $a\psi = XBV$.  To see that $\psi$ is well
defined, first observe that
\begin{equation}\label{stabilizing}
XA=X,\ YC=Y,\ AU=U,\ CV=V.
\end{equation}
  Using this we calculate
\begin{equation}\label{kquotient}
\begin{pmatrix} X & Z\\
                  0 & Y\\\end{pmatrix}\begin{pmatrix} A & B \\
                                                 0 &
C\end{pmatrix}\begin{pmatrix} U & W \\
                                                 0 &
V\end{pmatrix} = \begin{pmatrix} XU & XW+XBV+ZV \\
                                                 0 &
YV\end{pmatrix}.
\end{equation}
Subtracting $XW+ZV$ (which doesn't depend on the choice of a
representative of $a$) from the upper right hand corner shows that
$a\psi$ is well defined.  In fact it is evident from
\eqref{kquotient} that $\psi$ is injective.  We show that $\psi$
is a morphism to the additive group of $M_{m,r}(K)$.  It clearly
sends the identity matrix to $0$. Also
if $a,b\in S$ with respective middle coordinates $$\begin{pmatrix} A & B \\
                                                 0 &
C\end{pmatrix}, \begin{pmatrix} A' & B' \\
                                                 0 &
C'\end{pmatrix}$$ then $a\psi + b\psi = XBV + XB'V.$ But
the middle coordinate of $ab$ is $$\begin{pmatrix} AA' & AB'+BC'\\
0  & CC'\end{pmatrix}.$$ So $(ab)\psi = X(AB'+BC')V = XAB'V+XBC'V
= XB'V+XBV$ since $a,b\in S$ (cf.\ \eqref{stabilizing}).  Hence
$S$ is isomorphic to a finite subgroup of the additive group of
$M_{m,r}(K)$. In particular, if $K$ is a field and $\mathrm{char}\
K=0$, then $S$ must be trivial; if $\mathrm{char}\ K=p$, then
$S\in \pv {Ab}(p)$. The lemma follows.
\end{proof}

\subsubsection*{Proof of Theorem~\ref{Thm:dbsd}}
Let $K$ be a field and \pv V be a variety.  Let \pv U be the
smallest variety containing $\pv V$ such $\ell \pv I\dsd \pv U=\pv
U$ if $\mathrm{char}\ K=0$ or $\pv {Ab}(p)\dsd \pv U=\pv U$ if
$\mathrm{char}\ K=p$. Proposition~\ref{lVimpliesLV} and
\eqref{malcevclosed} immediately implies $\pv U\subseteq \pv
{LG}_{K}\malce \pv V$. To prove the converse, we need the
following.

\begin{Lemma}
Suppose $M$ is a finite submonoid of $M_n(K)$ in block upper
triangular form with diagonal block monoids $M_1,\ldots, M_k$
belonging to \pv V.  Then $M\in \pv U$.
\end{Lemma}
\begin{proof}
We induct on $k$.  If $k=1$, then $M=M_1\in \pv V\subseteq \pv U$.
In general, note that we can repartition $n$ into two blocks, one
corresponding to the union of the first $k-1$ of our original blocks
and the other corresponding to the last block.  We then obtain a
block upper triangular matrix monoid with two diagonal block
monoids $M'$ and $M_k$. By induction $M'\in \pv U$ (being block
upper triangular with $k-1$ diagonal blocks $M_1,\ldots, M_{k-1}$
belonging to \pv V) whilst $M_k\in \pv V\subseteq \pv U$.
Therefore $M'\times M_k\in \pv U$. Lemma~\ref{lemma:decomp} shows
that the kernel category of the projection to $M'\times M_k$
belongs to $\ell \pv I$, respectively $\ell \pv {Ab}(p)$,
according to the characteristic of $K$.  Hence $M\in \pv U$.
\end{proof}

To complete the proof of Theorem~\ref{Thm:dbsd}, suppose $M\in \pv
{LG}_{K}\malce \pv V$.  Consider the regular representation of $M$.
By finding a composition series for $M$, we can put $M$ in block
upper triangular form where the diagonal blocks $M_1,\ldots, M_k$
are the action monoids of the irreducible representations of $M$
over $K$. Since, by Theorem~\ref{malcevbymorph}, $M/\Radf K M\in \pv
V$, the $M_i$ belong to \pv V. The previous lemma then shows that
$M\in \pv U$, establishing Theorem~\ref{Thm:dbsd}.

\section{Applications to Formal Language Theory}
Another application of the Rhodes radical is to Formal Language
Theory, namely to unambiguous marked products and marked products
with counter.  Some of these results were announced in \cite{words}.

Recall that a word $u$ over a finite alphabet $\Sigma$ is said to
be a \emph{subword} of a word $v\in\Sigma^*$ if, for some $n\ge
1$, there exist words $u_1,\dots,u_n,v_0,v_1,\dots,v_n\in\Sigma^*$
such that $u=u_1u_2\cdots u_n$ and
\begin{equation}
\label{subword relation} v=v_0u_1v_1u_2v_2\cdots u_nv_n.
\end{equation}

The subword relation reveals interesting combinatorial properties
and plays a prominent role in formal language theory, as well as in
the theory of Coxeter groups via its relation to the Bruhat order
\cite{Coxeter}. For instance, recall that languages consisting of
all words over $\Sigma$ having a given word $u\in\Sigma^*$ as a
subword serve as a generating system for the Boolean algebra of
so-called \emph{piecewise testable} languages. It was a deep study
of combinatorics of the subword relation that led
Simon~\cite{Si72,Si75} to his elegant algebraic characterization of
piecewise testable languages. Further, the natural idea to put
certain rational constraints on the factors $v_0,v_1,\dots,v_n$ that
may appear in a decomposition of the form~\eqref{subword relation}
gave rise to the useful notion of a marked product of languages
studied from the algebraic viewpoint by
Sch\"utzenberger~\cite{Sch65}, Reutenauer~\cite{Reutenauer},
Straubing~\cite{Straubing}, Simon~\cite{Simonpaper}, amongst others.

Yet another natural idea is to count how many times a word
$v\in\Sigma^*$ contains a given word $u$ as a subword, that is, to
count different decompositions of the form~\eqref{subword
relation}. Clearly, if one wants to stay within the realm of
rational languages, one can only count up to a certain threshold
and/or modulo a certain number. For instance, one may consider
Boolean combinations of languages consisting of all words over
$\Sigma$ having $t$ modulo $p$ occurrences of a given word
$u\in\Sigma^*$ (where $p$ is a given prime number). This class of
languages also admits a nice algebraic characterization,
see~\cite[Sections VIII.9 and VIII.10]{Eilenberg} and
also~\cite{Th83}. Combining modular counting with rational
constraints led to the idea of marked products with modular
counters explored, in particular, by Weil~\cite{Weil} and
P\'eladeau~\cite{collapse}.

The most natural version of threshold counting is formalized via
the notion of an unambiguous marked product in which one considers
words $v\in\Sigma^*$ having exactly one
decomposition~\eqref{subword relation} with a given subword $u$
and given rational constraints on the factors $v_0,v_1,\dots,v_n$.
Such unambiguous marked products have been investigated by
Sch\"utzenberger~\cite{Sch76}, Pin~\cite{Pinunambi}, Pin,
Straubing, and Th\'erien~\cite{Pinetal}, amongst others.

Many known facts on marked products rely on rather difficult
techniques from finite semigroup theory, namely, on the bilateral
semidirect product decomposition results of Rhodes \emph{et al.}
\cite{Kernel,MPS2} mentioned above. These results are proved
using Rhodes's classification of maximal proper
surmorphisms~\cite{MPS,MPS2,Arbib} via case-by-case analysis of the
kernel categories of such maps~\cite{Kernel,MPS2}. The aim of
the present section is to give easier and --- we hope --- more
conceptual proofs of several crucial facts about marked products by
using matrix representations of finite semigroups as a main tool. In
particular, we are able to prove the results of P\'eladeau and Weil in
one step, without any case-by-case analysis and without using the
machinery of categories. Rather we adapt Simon's analysis of the
combinatorics of multiplying upper triangular
matrices~\cite{Simonpaper} from the case of Sch\"utzenberger
products to block upper triangular matrices.  We failed to obtain
such a purely combinatorial argument for the case of unambiguous
products; we still need to use a lemma on kernel categories.
Nevertheless we have succeeded in avoiding the decomposition results
and case-by-case analysis.

Recall that Eilenberg established~\cite[Vol.~B,
Chap.~VII]{Eilenberg} a correspondence between varieties of finite
monoids and so-called varieties of languages. If \pv V is a variety
of finite monoids and $\Sigma$ a finite alphabet, then ${\mathcal
V}(\Sigma^*)$ denotes the set of all languages over $\Sigma$ that
can be recognized by monoids in \pv V. (Such languages are often
referred to as \pv V-\emph{languages}.) The operator $\mathcal V$
that assigns each free monoid $\Sigma^*$ the set ${\mathcal
V}(\Sigma^*)$ is said to be \emph{the variety of languages
associated to} \pv V. The syntactic monoid~\cite[loc.\
cit.]{Eilenberg} of a rational language $L$ will be denoted $M_L$.
It is known that $L$ is a \pv V-language if and only if $M_L\in \pv
V$.

\subsection{Products with Counter}
Our first application is to prove the results of P\'eladeau and
Weil~\cite{collapse,Weil} on products with counter.

Let $L_0,\ldots, L_m\subseteq\Sigma^*$, $a_1,\ldots,a_m\in\Sigma$
and let $n$ be an integer.  Then the \emph{marked product with
modulo $n$ counter} $L=(L_0a_1L_1\cdots a_mL_m)_{r,n}$ is the
language of all words $w\in \Sigma^*$ with $r$ factorizations
modulo $n$ of the form $w=u_0a_1u_1\cdots a_mu_m$ with each
$u_i\in L_i$. One can show that $L$ is rational~\cite{Weil} (see
also the proof of Theorem~\ref{Pel-Weil} below). Using a
decomposition result of Rhodes and Tilson~\cite{Kernel} (see also
\cite{MPS2}) based on case-by-case analysis of kernel
categories of maximal proper surmorphisms (see
\cite{MPS,MPS2,Arbib}), Weil characterized the closure of a
variety $\mathcal V$ under marked products with modulo $p$
counter.  This required iterated usage of the so-called ``block
product'' principle. But Weil missed that the Boolean algebra
generated by $\mathcal V(\Sigma^*)$ and marked products with
modulo $p$ counters of members $\mathcal V(\Sigma^*)$ is already
closed under marked products with modulo $p$ counters; this was
later observed by P\'eladeau~\cite{collapse}. The difficulty arises
because it is not so clear how to combine marked products with
modulo $p$ counters into new marked products with modulo $p$
counters.

We use representation theory to prove the result in one fell
swoop.   Our approach is inspired by a paper of
Simon~\cite{Simonpaper} dealing with marked products and the
Sch\"utzenberger product of finite semigroups.

\begin{Lemma}\label{blockupcharp}
Let \pv V be a variety of finite monoids, $\p:\Sigma^*\to M$ be a
morphism with $M$ finite.  Let $K$ be a field of characteristic
$p$ and suppose that $M$ can be represented faithfully by block
upper triangular matrices over $K$ so that the monoids formed by
the diagonal blocks of the matrices in the image of $M$ all belong
to \pv V. Let $F\subseteq M$. Then  $L=F\pinv$ is a Boolean
combination of members of $\mathcal V(\Sigma^*)$ and of marked
products with modulo $p$ counter $(L_0a_1L_1\cdots a_nL_n)_{r,p}$
with the $L_i\in \mathcal V(\Sigma^*)$.
\end{Lemma}

\begin{proof}
Suppose $M\leq M_t(K)$ and $t=t_1+\cdots +t_k$ is the partition of
$t$ giving rise to the block upper triangular form. Let $M_i$ be
the monoid formed by the $t_i\times t_i$ matrices over $K$ arising
as the $i^{th}$ diagonal blocks of the matrices in the image of
$M$. Given $w\in \Sigma^*$ and $i,j\in \{1,\ldots,k\}$, define
$\p_{i,j}:\Sigma^*\to M_{t_i,t_j}(K)$ by setting $w\p_{i,j}$ to be
the $t_i\times t_j$ matrix that is the $i,j$-block of the block
upper triangular form. So in particular $w\p_{i,j}=0$ for $j<i$.
Also $\p_{i,i}$ is a morphism $\p_{i,i}:\Sigma^*\to M_i$ for all
$i$.

First we observe that we may take $F$ to be a singleton $\{u\p\}$.
For each $1\leq i\leq j\leq k$, let $$L_{i,j}=\{w\in \Sigma^*\mid
w\p_{i,j} = u\p_{i,j}\}.$$  Then clearly
$$u\p\pinv = \bigcap_{1\leq i\leq j\leq k} L_{i,j}.$$
Since $L_{i,i}$ is recognized by $M_i$, it suffices to show
$L_{i,j}$, where $1\leq i<j\leq k$, can be written as a Boolean
combination of marked products with modulo $p$ counter of
languages recognized by the $M_l$. Changing notation, it suffices
to show that if $1\leq i<j\leq k$ and $C\in M_{t_i,t_j}(K)$, then
\begin{equation}\label{langC}
L(C)=\{w\in \Sigma^*\mid w\p_{i,j}=C\}
\end{equation}
is a Boolean combination of marked products with modulo $p$
counter of languages recognized by the $M_i$.

The following definitions are inspired by~\cite{Simonpaper},
though what Simon terms an ``object'', we term a ``walk''.  A
\emph{walk} from $i$ to $j$ is a sequence
\begin{equation}\label{walk}
\mathfrak w =(i_0,m_0,a_1,i_1,m_1,\dots,a_r,i_r,m_r)
\end{equation}
where $i=i_0<i_1<\cdots<i_r=j$, $a_l\in \Sigma$ and $m_l\in
M_{i_l}$. There are only finitely many walks.  The set of walks
will be denoted $\mathfrak W$. Given a walk $\mathfrak w$, we
define its \emph{value} to be
$$\mathsf{v}(\mathfrak w) =m_0(a_1\p_{i_0,i_1})m_1\cdots
(a_r\p_{i_{r-1},i_r})m_r\in M_{t_i,t_j}(K).$$ If $\mathfrak w$ is
a walk, we define the \emph{language} of $\mathfrak w$ to be the
marked product
$$L(\mathfrak w) = (m_0\p_{i_0,i_0}\inv) a_1(m_1\p_{i_1,i_1}\inv)
\cdots a_r(m_r\p_{i_r,i_r}\inv).$$

If $w\in \Sigma^*$ and $\mathfrak w$ is a walk of the form
\eqref{walk}, we define $w(\mathfrak w)$ to be the
\emph{multiplicity} of $w$ in $L(\mathfrak w)$, that is, the
number of factorizations $w=u_0a_1u_1\cdots a_ru_r$ with
$u_l\p_{i_l,i_l} = m_l$; this number is taken to be $0$ if there
are no such factorizations.  If $0\leq n<p$, we establish the
shorthand
$$L(\mathfrak w)_{n,p} =\left((m_0\p_{i_0,i_0}\inv) a_1(m_1\p_{i_1,i_1}\inv)\cdots
a_r(m_r\p_{i_r,i_r}\inv)\right)_{n,p}.$$ Notice  that $L(\mathfrak
w)_{n,p}$ consists of all words $w$ with $w(\mathfrak w)\equiv
n\bmod p$ and is a marked product with modulo $p$ counter of
$\mathcal V(\Sigma^*)$ languages.

The following is a variant of~\cite[Lemma~7]{Simonpaper}.

\begin{Claim}\label{factorizing}
Let $w\in \Sigma^*$.  Then
\begin{equation}\label{walkeq}
w\p_{i,j} = \sum_{\mathfrak w\in \mathfrak W} w(\mathfrak
w)\mathsf{v}(\mathfrak w).
\end{equation}
\end{Claim}

\begin{proof}
Let $w=b_1\cdots b_r$ be the factorization of $w$ in letters. Then
the formula for matrix multiplication gives
\begin{equation}\label{matrixmult}
w\p_{i,j} = \sum
(b_1\p_{i_0,i_1})(b_2\p_{i_1,i_2})\cdots(b_r\p_{i_{r-1},i_r})
\end{equation}
where the sum extends over all $i_l$ such that $i_0=i$, $i_r=j$
and $i_l\in\{1,\ldots,k\}$ for $0<l<r$. Since $v\p_{l,n} =0$ for
$l>n$, it suffices to consider sequences such that $i=i_0\leq
i_1\leq \cdots \leq i_r=j$.   For such a sequence, we may group
together neighboring indices that are equal.  Then since all the
$\p_{n,n}$ are morphisms, we see that each summand in
\eqref{matrixmult} is the value of a walk $\mathfrak w$ and that
$\mathfrak w$ appears exactly $w(\mathfrak w)$ times in the sum.
\end{proof}

To complete the proof of Lemma~\ref{blockupcharp}, we observe that $L(C)$ (defined in
\eqref{langC}) is a Boolean combination of languages of the form
$L(\mathfrak w)_{n,p}$.  Let $X$ be the set of all functions
$f:\mathfrak W\to \{0,\ldots,p-1\}$ such that $$\sum _{\mathfrak
w\in \mathfrak W} f(\mathfrak w)\mathsf v(\mathfrak w) = C.$$ It
is then immediate from \eqref{walkeq} and $\mathrm{char} K=p$ that
$$L(C) = \bigcup_{f\in X} \bigcap_{\mathfrak w\in\mathfrak W}
L(\mathfrak w)_{f(\mathfrak w),p}$$ completing the proof.
\end{proof}

\begin{Thm}\label{Pel-Weil}
Let $L\subseteq \Sigma^*$ be a rational language, \pv V be a
variety of finite monoids and $K$ be a field of characteristic
$p$. Then the following are equivalent.
\begin{enumerate}
\item $M_L\in \pv{LG}_p\malce \pv V$; \item $M_L/\Radf {K}
{M_L}\in \pv V$; \item $M_L$ can be faithfully represented by
block upper triangular matrices over $K$ so that the monoids
formed by the diagonal blocks of the matrices in the image of
$M_L$ all belong to \pv V; \item $L$ is a Boolean combination of
members of $\mathcal V(\Sigma^*)$ and languages $(L_0a_1L_1\cdots
a_nL_n)_{r,p}$ with the $L_i\in \mathcal V(\Sigma^*)$.
\end{enumerate}

\end{Thm}
\begin{proof}
The equivalence of (1) and (2) was established in Theorem
\ref{malcevbymorph}.

For (2) implies (3), take a composition series for the regular
representation of $M_L$ over $K$: it is then in block upper
triangular form and, by (2) the monoids formed by diagonal blocks
of matrices in the image of $M_L$ all belong to \pv V, being the
action monoids from the irreducible representations of $M_L$ over
$K$.

(3) implies (4) is immediate from Lemma~\ref{blockupcharp}.

For (4) implies (1), it suffices to deal with a marked product with
counter $L=(L_0a_1L_1\cdots a_nL_n)_{r,p}.$ Let $\mathcal A_i$ be
the minimal trim deterministic automaton \cite[Vol.~A]{Eilenberg}
of $L_i$. Let $\mathcal A$ be the non-deterministic automaton
obtained from the disjoint union of the $\mathcal A_i$ by attaching
an edge labelled $a_i$ from each final state of $\mathcal A_{i-1}$
to the initial state of $\mathcal A_i$.  To each letter
$a\in\Sigma$, we associate the matrix $a\p$ of the relation that $a$
induces on the states. Since $a\p$ is a $0,1$-matrix, we can view it
as a matrix over $\FF_p$. In this way we obtain a morphism
$\p:\Sigma^*\to M_k(\FF_p)$ where $k$ is the number of states of
$\mathcal A$. Let $M = \Sigma^*\p$. Trivially, $M$ is finite. We
observe that $M$ is  block upper triangular with diagonal blocks the
syntactic monoids $M_{L_i}$ (the partition of $k$ arises from taking
the states of each $\mathcal A_i$). Notice that $M$ recognizes $L$,
since $L$ consists of all words $w$ such that $(w\p)_{s,f}=r$ where
$s$ is the start state of $\mathcal A_0$ and $f$ is a final state of
$\mathcal A_n$. Applying Lemma~\ref{toalgebra} to the projection to
the diagonal blocks gives that $M$ and its quotient $M_L$ belong to
$\pv{LG}_p\malce\pv V$.
\end{proof}

The proof of (4) implies (1) gives a fairly easy argument that
marked products of rational languages with mod $p$ counter are
rational.

Since the operator $\pv {LG}_p\malce (\ )$ is idempotent, we
immediately obtain the following result of~\cite{collapse,Weil}.

\begin{Cor}
Let \pv V be a variety of finite monoids and $\pv W=\pv
{LG}_p\malce \pv V$.  Let $\mathcal{W}$ be the corresponding
variety of languages. Then
\begin{enumerate}
\item $\mathcal{W}(\Sigma^*)$ is the smallest class of languages
containing $\mathcal {V}(\Sigma^*)$, which is closed under Boolean
operations and formation of marked products with modulo $p$
counters.

\item  $\mathcal{W}(\Sigma^*)$ consists of all Boolean
combinations of elements of $\mathcal V(\Sigma^*)$ and  marked
products with modulo $p$ counters of elements of $\mathcal
{V}(\Sigma^*)$.
\end{enumerate}
\end{Cor}

Some special cases are the following.  If $\pv V$ is the trivial
variety of monoids, then $\pv {LG}_p\malce \pv V=\pv {G}_p$ and we
obtain Eilenberg's result~\cite[Section~VIII.10]{Eilenberg} that
the $\pv {G}_p$ languages consist of the Boolean combinations of
languages of the form $(\Sigma^*a_1\Sigma^*\cdots
a_n\Sigma^*)_{r,p}$. Notice that $\pv {G}_p$ consists of the
groups unitriangularizable over characteristic $p$.  The languages
over $\Sigma^*$ associated to $\pv {LG}_p\malce \pv {Sl}$ (as
observed in~\cite{amv} and Theorem~\ref{Thm:unitriangularizability},
this variety consists of the
unitriangularizable monoids over characteristic $p$) are the
Boolean combinations of languages of the forms
$$\Sigma^*a\Sigma^*\ \text{ and }\ (\Sigma_0^*a_1\Sigma_1^*\cdots a_n\Sigma_n^*)_{r,p}$$
where $\Sigma_i\subseteq\Sigma$.

We remark that Weil shows~\cite{Weil} that closing $\mathcal
V(\Sigma^*)$ under marked products with modulo $p^n$ counters, for
$n>1$, does not take you out of the $\pv {LG}_p\malce \pv
V$-languages.

\subsection{Unambiguous Products}

Our next application is to recover results of Sch\"utzenberger, Pin,
Straubing, and Th\'erien concerning unambiguous products. Our proof
of one direction is along the lines of~\cite{Pinetal} but our usage
of representation theory allows us to avoid using results relying on
case-by-case analysis of maximal proper surmorphisms.

Let $\Sigma$ be a finite alphabet, $L_0,\dots,L_n\subseteq
\Sigma^*$ be rational languages and $a_1,\ldots,a_n\in\Sigma$.
Then the \emph{marked product} $L=L_0a_1L_1\cdots a_nL_n$ is
called \emph{unambiguous} if each word $w\in L$ has exactly one
factorization of the form $u_0a_1u_1\cdots a_nu_n$, where each
$u_i\in L_i$.  We also allow the degenerate case $n=0$.

We shall need to use a well-known and straightforward consequence
of the distributivity of concatenation over union
(cf.~\cite{Pinetal}), namely, if $L_0,\ldots, L_n$ are
disjoint unions of unambiguous marked products of elements of
$\mathcal V(\Sigma^*)$, then the same is true for any unambiguous
product $L_0a_1L_1\cdots a_nL_n$.  We also need a lemma about
languages recognized by finite monoids of block upper triangular
matrices in characteristic $0$.

\begin{Lemma}\label{blockupchar0}
Let \pv V be a variety of finite monoids, $\p:\Sigma^*\to M$ be a
morphism with $M$ finite.  Let $K$ be a field of characteristic
$0$ and suppose that $M$ can be represented faithfully by block
upper triangular matrices over $K$ so that the monoids
$M_1,\dots,M_k$ formed by diagonal blocks of matrices in the image
of $M$ all belong to \pv V. Let $F\subseteq M$. Then $L=F\pinv$ is
a disjoint union of unambiguous marked products $L_0a_1L_1\cdots
a_nL_n$ with the $L_i\in \mathcal V(\Sigma^*)$.
\end{Lemma}

\begin{proof}
We induct on the number $k$ of diagonal blocks. If there is only
one block we are done.

Now let $k>1$. We can repartition $n$ into two blocks, one
corresponding to the union of the first $k-1$ of our original
blocks and the other corresponding to the last block. The first
diagonal block, call it $N$, is block upper triangular with
diagonal blocks $M_1,\dots, M_{k-1}$; the second is just $M_k$. By
induction, any language recognized by $N$ is a disjoint union of
unambiguous marked products $L_0a_1L_1\cdots a_rL_r$ with the
$L_i\in \mathcal V(\Sigma^*)$.  Since $M_k\in \pv V$ it is easy to
check that any language recognized by $N\times M_k$ is also a
disjoint union of unambiguous marked products $L_0a_1L_1\cdots
a_rL_r$ with the $L_i\in \mathcal V(\Sigma^*)$. Thus to prove the
result, it suffices to show that $L$ is a disjoint union of
unambiguous marked products $L_0a_1L_1\cdots a_nL_n$ with the
$L_i$ recognized by $N\times M_k$. By Lemma~\ref{lemma:decomp},
the projection from $M$ to $N\times M_k$ has locally trivial
kernel category. Then
\cite[Proposition 2.2]{Pinetal} shows us that $L$ is a disjoint
union of such unambiguous marked products.
\end{proof}

We ask whether there is a simple combinatorial proof of this lemma
that avoids the use of~\cite[Proposition 2.2]{Pinetal} along the
lines of the proof of Lemma~\ref{blockupcharp}.

\begin{Thm}
Let $L\subseteq \Sigma^*$ be a rational language, \pv V be a
variety of finite monoids and $K$ a field of characteristic $0$.
Then the following are equivalent.
\begin{enumerate}
\item $M_L\in \pv{LI}\malce \pv V$; \item $M_L/\Radf {K} {M_L}\in
\pv V$; \item $M_L$ can be faithfully represented by block upper
triangular matrices over $K$ so that the monoids formed by the
diagonal blocks of the matrices in the image of $M_L$ all belong
to \pv V.\; \item $L$ is a disjoint union of unambiguous products
$L_0a_1L_1\cdots a_nL_n$ with the $L_i\in \mathcal V(\Sigma^*)$.
\end{enumerate}
\end{Thm}
\begin{proof}
The equivalence of (1) and (2) follows from Theorem
\ref{malcevbymorph}.

For (2) implies (3), take a composition series for the regular
representation of $M_L$ over $K$: it is then in block upper
triangular form and by (2) monoids formed by diagonal blocks of
matrices in the image of $M_L$ all belong to \pv V.

(3) implies (4) is immediate from Lemma~\ref{blockupchar0}.

For (4) implies (1), it suffices to deal with a single unambiguous
marked product $L=L_0a_1L_1\cdots a_nL_n$. Let $\mathcal A_i$ be the
minimal trim deterministic automaton for $L_i$ and let $\mathcal A$
be the non-deterministic automaton obtained from the disjoint union
of the $L_i$ by attaching an edge labelled $a_i$ from each final
state of $\mathcal A_{i-1}$ to the initial state of $\mathcal A_i$.
To each letter $a\in A$, we associate the matrix $a\p$ of the
relation that $a$ induces on the states. In this way we obtain a
morphism $\p:\Sigma^*\to M_k(\QQ)$ where $k$ is the number of states
of $\mathcal A$. Let $M = \Sigma^*\p$. We observe that $M$ is block
upper triangular with diagonal blocks the syntactic monoids
$M_{L_i}$ (the partition of $k$ arises from taking the states of
each $\mathcal A_i$). Notice that $M$ recognizes $L$, since $L$
consists of all words $w$ such that $(w\p)_{s,f}>0$ where $s$ is the
start state of $\mathcal A_0$ and $f$ is a final state of $\mathcal
A_n$.  First we show that $M$ is finite.  In fact, we claim $M$
contains only $0,1$-matrices (and hence must be finite). Indeed,
suppose $(w\p)_{i,j}>1$ some $i,j$. Since each $M_{L_i}$ consists of
$0,1$-matrices, we must have that $i$ is a state of some $\mathcal
A_l$ and $j$ a state of some $\mathcal A_r$ with $l<r$. But
$(w\p)_{i,j}$ is the number of paths labelled by $w$ from $i$ to $j$
in $\mathcal A$.  Thus if $u,v$ are words reading respectively from
the start state of $\mathcal A_0$ to $i$ and from $j$ to a final
state of $\mathcal A_n$ (such exist since the $\mathcal A_i$ are
trim), then $uwv$ has at least two factorizations witnessing
membership in $L$, contradicting that $L$ was unambiguous.  Since
the collection of all block upper triangular matrices is an algebra
over $\QQ$, as is the collection of block diagonal matrices, an
application of Lemma~\ref{toalgebra} to the projection to the
diagonal blocks gives that $M\in \pv {LI}\malce \pv V$ and so, since
$M\onto M_L$, we have $M_L\in \pv {LI}\malce \pv V$.
\end{proof}

Since the operator $\pv {LI}\malce (\ )$ is idempotent, we
immediately obtain the following result
of~\cite{Pinunambi,Pinetal}.
\begin{Cor}
Let \pv V be a variety of finite monoids and $\pv W=\pv {LI}\malce
\pv V$.  Let $\mathcal{W}$ be the corresponding variety of
languages. Then
\begin{enumerate}
\item $\mathcal{W}(\Sigma^*)$ is the smallest class of languages
containing $\mathcal {V}(\Sigma^*)$, which is closed under Boolean
operations and formation of unambiguous marked products. \item
$\mathcal{W}(\Sigma^*)$ consists of all finite disjoint unions of
unambiguous marked products of elements of $\mathcal
{V}(\Sigma^*)$.
\end{enumerate}
\end{Cor}

Recall that the Malcev product of the pseudovariety \pv{LI} with
the pseudovariety $\pv{Sl}$ of semilattices
(idempotent commutative monoids) is equal to the famous
pseudovariety \pv{DA} of all finite monoids whose regular
$\mathcal{D}$-classes are idempotent subsemigroups (see~\cite{TT}
for a nice survey of combinatorial, logical and automata-theoretic
characterizations of  \pv{DA}). Applying the above corollary, one
obtains the classical result of Sch\"utzenberger~\cite{Sch76} that
$\mathcal{DA}(\Sigma^*)$ consists of disjoint unions of
unambiguous products of the form $\Sigma_0^*a_1\Sigma_1^*\cdots
a_n\Sigma_n^*$ with $\Sigma_i\subseteq \Sigma$ for all $i$. We saw
in Corollary~\ref{DAisforever} that $\pv {DA}$ consists of
precisely those finite monoids that can be faithfully represented
by upper triangular matrices with zeroes and ones on the diagonal
over $\mathbb{Q}$.

\section{\v{C}ern\'y's Conjecture for \pv{DS}}
A deterministic automaton $\mathcal A = (Q,A)$ is called
\emph{synchronizing} if there is a word $w\in A^*$ such that
$|Qw|=1$, that is $w$ acts as a constant map on $Q$.  Such a word
$w$ is called a \emph{synchronizing word}.  \v{C}ern\'y raised the
following question: how large can a minimal length synchronizing
word for a synchronizing automaton be as a function of the number of
states of the automaton?  He showed that for each $n>1$, there are
$n$ state synchronizing automata with minimal synchronizing words of
size $(n-1)^2$~\cite{Cerny}.  The best known upper bound, due to Pin
\cite{pincerny}, is $\frac{n^3-n}{6}$.  \v{C}ern\'y conjectured
that in fact $(n-1)^2$ is the exact answer. Many special cases of
the conjecture have been proved (for instance, \cite{pincirc,Dubuc,Kari, Volkov}),
but the conjecture in general remains wide open.

In this section we show, using representation theory, that
\v{C}ern\'y's conjecture is true for synchronizing automata with
transition monoids in the variety \pv {DS}.  We begin by giving a
representation theoretic rephrasing of the problem from the thesis
of Steinberg's Master's student Arnold~\cite{rick}.

Let $\mathcal A=(Q,A)$ be a deterministic automaton and let $M$ be
its transition monoid.  Set $n=|Q|$.   Let $V$ be the $\mathbb
Q$-vector space with basis \mbox{$B=\{e_q\mid q\in Q\}$}. Then there is a
faithful representation $\p:M\to \End {\mathbb Q} V$ defined on
the basis by $$e_qm\p = e_{qm}.$$  We consider $V$ with the usual
inner product.  Let
$$V_1 = \mathrm{Span}\{\sum_{q\in Q} e_q\}\ \text{and}\
V_0=V_1^{\perp}.$$ We claim that $V_0$ is $M$-invariant.  Indeed,
suppose $v\in V_0$ and $m\in M$.  Let $v_1=\sum_{q\in Q} e_q$.
Then $$\langle vm\p,v_1\rangle = \langle v,v_1(m\p)^T\rangle$$
(where $()^T$ denotes transposition).  With respect to the basis $B$,
$m\p$ is a row monomial matrix (meaning each row has precisely one
non-zero entry) and hence $m\p^T$ is column monomial.  On the
other hand, in the basis $B$, the vector $v_1$ is the vector of
all ones, hence $v_1$ is fixed by any column monomial matrix. Thus
$$\langle vm\p,v_1\rangle = \langle v,v_1(m\p)^T\rangle = \langle
v,v_1\rangle=0,$$ establishing that $vm\p \in V_0$, as desired.
We conclude that $V_0$ is $M$-invariant.  Let $\psi:M\to \End
{\QQ} {V_0}$ be the associated representation.

Without loss of generality, let us assume $Q=\{1,\ldots, n\}$.
Then $V_0$ has basis $B_0=\{f_1,\ldots, f_{n-1}\}$ where $f_i =
e_n-e_i$.  Also $f_im\psi = e_{n\cdot m}-e_{i\cdot m}$.  In
particular, $m\psi =0$ if and only if $n\cdot m = i\cdot m$ for
all $i=1,\ldots, n-1$, that is, if and only if $m$ is a constant
map.  Thus $w\in A^*$ is a synchronizing word if and only if
$[w]\psi =0$, where $[w]$ is the image of $w$ in $M$.  In
particular, $\mathcal A$ is synchronizing if and only if $M\psi$
contains the zero matrix. Since $V_0$ has dimension $n-1$, we will
have proved \u{C}ern\'y's conjecture for the case that $M$ belongs to \pv
{DS} once we have proven the following theorem, which can be viewed as
the ``mortality problem'' for \pv {DS}.

\begin{Thm}\label{DSmortality}
Let $K$ be a field and let $A$ be a finite alphabet.  Let $M$ be a
finite $A$-generated submonoid of $M_k(K)$ belonging to \pv {DS}
and suppose that $0\in M$. Then there exists a word $w\in A^*$ of
length at most $k^2$ such that $w$ maps to $0$ in $M$.
\end{Thm}

Before proving this theorem, we need a lemma.

\begin{Lemma}\label{DSGGM}
Let $S\in \pv {DS}$ be a non-trivial generalized group mapping
semigroup with a zero element $0$.  Then $S\setminus \{0\}$ is a
subsemigroup.
\end{Lemma}
\begin{proof}
By definition, $S$ has a ($0$-)minimal ideal $I$ on which it acts
faithfully on both the left and right.  Since $S$ is non-trivial,
$I$ cannot be the ideal $0$.  Thus $I$ is $0$-minimal.  Since $I$
is regular, $I\setminus \{0\}$ is a regular $\J$-class $J$.
Suppose $s,t\geq_{\J} J$.  Then, since $S\in \pv {DS}$, we have
$st\geq _{\J} J$, see~\cite[Section~8.1]{almbook}.  Since $S$ acts
faithfully on $I$, only $0$ is not $\J$-above $J$.  Thus
$S\setminus \{0\}$ is indeed a subsemigroup.
\end{proof}

\subsubsection*{Proof of Theorem~\ref{DSmortality}}
By choosing a composition series for the $KM$-module $K^k$, we can
place $M$ in block upper triangular form where the diagonal block
monoids $M_1,\ldots, M_r$, with $1\leq r\leq k$, are irreducible.
Since each $M_i$ is a homomorphic image of $M$, they each have a
zero element and each belong to \pv {DS}.  Being irreducible, they
are generalized group mapping monoids by Theorem
\ref{thm:irredGGM}.  Thus $M_i\setminus \{0\}$ is a submonoid by
Lemma~\ref{DSGGM}.   Let $\alpha:M\to M_1\times\cdots\times M_r$
be the projection.  Suppose $w\in A^*$ maps to zero in $M$, then
$w\alpha =0$ and hence, for each $i=1,\ldots, r$, there is a
letter $a_i\in A$ with the $i^{th}$ coordinate of $a_i\alpha$
equal to zero (using that the product of non-zero elements of
$M_i$ remains non-zero).  Thus we can find a word $u\in A^*$ of
length at most $r\leq k$ such that $u$ represents an element $m$
of $M$ with zeroes on the diagonal blocks.  But then $m$ is
nilpotent of index at most $k$ since it is a $k\times k$ upper
triangular matrix with zeroes on the diagonal.  Thus $u^k$
represents $0$ and $|u^k| \leq k^2$.\qed

\medskip

We remark that the proof gives a bound of $\min \{|A|,r\}\cdot r$
where $r$ is the number of irreducible constituents of $M$.  This is
because in forming $u$ we do not need to repeat letters and because
the nilpotency index is actually bounded by the number of zero
blocks on the diagonal. Hence if either $|A|$ or $r$ are small, then
we can do better than $k^2$.

Applying the above theorem in the context of the representation
$\psi$ of the transition monoid of an automaton on $V_0$ discussed
above, we obtain the following theorem, verifying \u{C}ern\'y's
conjecture for \pv {DS}.

\begin{Thm}
Every synchronizing automaton on $n$ states with transition monoid
in \pv {DS} has a synchronizing word of length at most $(n-1)^2$.
\end{Thm}

We do not know whether $(n-1)^2$ is sharp when restricted to
automata with transition monoids in \pv {DS}.

A further application of the representation theory to \v{C}ern\'y's conjecture
can be found in a recent paper by  F.~Arnold and the third author~\cite{Arnold-Steinberg}.

\bibliographystyle{amsplain}

\end{document}